\documentclass[12pt]{article}
\usepackage{amssymb, epsfig, amsmath,amsfonts}

\usepackage{rotating}
\newlength{\defbaselineskip}
\newcommand{\setlinespacing}[1]%
           {\setlength{\baselineskip}{#1 \defbaselineskip}}

\newcommand{\singlespacing}{\setlength{\baselineskip}{\defbaselineskip}}

\RequirePackage[OT1]{fontenc}
\RequirePackage{amsthm,amsmath,epsfig,amssymb}
\usepackage{amsmath,amsfonts,amssymb}
\RequirePackage[numbers]{natbib}
\RequirePackage[colorlinks,citecolor=blue,urlcolor=blue]{hyperref}

\numberwithin{equation}{section}
\theoremstyle{plain}

\usepackage{amssymb, epsfig, amsmath,amsfonts}

\usepackage{rotating}

\newtheorem{lemma}{Lemma}[section]
\newtheorem{theorem}{Theorem}[section]
\newtheorem{definition}{Definition}[section]
\newtheorem{coro}{Corollary}[section]
\newtheorem{prop}{Proposition}[section]
\newtheorem{remark}{Remark}[section]

\newcommand{\bsq}{\vrule height .9ex width .8ex depth -.1ex}
\newcommand{\RR}{{\mathbb R}}

\newcommand{\sD}{{\cal D}}

\newcommand{\sC}{{\cal C}}

\newcommand{\beql}[1]{\begin{equation}\label{#1}}
\newcommand{\eeq}{\end{equation}}

\newcommand{\beqal}[1]{\begin{eqnarray}\label{#1}}
\newcommand{\eeqa}{\end{eqnarray}}
\newcommand{\beq}{\begin{displaymath}}
\newcommand{\eeqno}{\end{displaymath}}

\newcommand{\la}{\lambda}

\newcommand{\bone}{{\mathbf 1}}

\newcommand{\qandq}{\quad\mbox{and}\quad}
\newcommand{\qifq}{\quad\mbox{if}\quad}

\newcommand{\non}{\nonumber}

\begin{document}

%\renewcommand{\baselinestretch}{1.56}
%\setlength{\textwidth}{6.5in} \setlength{\textheight}{9in}
%\topmargiN=0.0in \headheight=0.0in \headsep=0.0in \oddsidemargin
%-0.1in
%\parskip 6pt

\title{\Large{Equivalence of Fluid Models for $G_{\lowercase{t}}/GI/N+GI$ Queues}}\vskip .1in
\author{
{\sf Weining Kang}\thanks{Department of Mathematics and Statistics, University of Maryland, Baltimore County, Baltimore, MD 21250 
 (\texttt{wkang@umbc.edu})}
~~~~
{\sf Guodong Pang}\thanks{The Harold and Inge Marcus Department of Industrial and Manufacturing Engineering,
Pennsylvania State University, University Park, PA 16802 (\texttt{gup3@psu.edu})} 
}

\date{\vspace{-5ex}}

\maketitle

%\begin{frontmatter}
%
%\title{Equivalence of Fluid Models for $G_{\lowercase{t}}/GI/N+GI$ Queues}
%\runtitle{Equivalence of Fluid Models for $G_{\lowercase{t}}/GI/N+GI$ Queues}
%
%%\begin{aug}
%\author{\fnms{Weining} \snm{Kang}\ead[label=e1]{wkang@umbc.edu}}
%\address{Department of Mathematics and Statistics, \\ University of Maryland, Baltimore County, \\
%Baltimore, MD 21250 \\
%\printead{e1}}
%and
%\author{\fnms{Guodong} \snm{Pang}
%\ead[label=e2]{gup3@psu.edu}}
%
%\address{The Harold and Inge Marcus Department  \\of Industrial and Manufacturing Engineering,\\
%Pennsylvania State University, \\
%University Park, PA 16802 \\
%\printead{e2}}
%
%\runauthor{ Kang and Pang}
%
%\affiliation{University of Maryland at Baltimore County \\ and The Pennsylvania State University}
%
%

\allowdisplaybreaks

%\doublespacing

\begin{abstract}

Four different fluid model formulations have been recently developed for $G_t/GI/N+GI$ queues, including a two-parameter fluid model in Whitt (2006) by tracking elapsed service and patience times of each customer, a measure-valued fluid model in Kang and Ramanan (2010) and its extension in Zu{\~n}iga (2014) by tracking elapsed service and patience times of each customer,  and a measure-valued fluid model in Zhang (2013) by tracking residual service and patience times of each customer. 
We show that the two fluid models tracking elapsed times (Whitt's and Kang and Ramanan's fluid models) are equivalent formulations for the same $G_t/GI/N+GI$ queue, whereas Zu{\~n}iga's fluid model and Zhang's fluid model are not entirely equivalent under general initial conditions. We then identify necessary and sufficient conditions under which Zu{\~n}iga's fluid model and Zhang's fluid model can be derived from each other for the same system, in which certain measure-valued fluid processes tracking residual service and patience times of each customer derived from Kang-Ramanan and Zu{\~n}iga's fluid models play an important role.
 The equivalence properties discovered provide important implications for the understanding of the recent development for non-Markovian many-server queues.

\end{abstract}

%This measure-valued fluid model tracking residual times is motivated by a 
%coupling property that connects it with fluid models tracking elapsed times. This coupling property enables us to construct a solution to the measure-valued fluid model tracking residual ({\it resp.} elapsed) times from a solution to the measure-valued fluid model tracking elapsed ({\it resp.} residual) times, and thus prove the well-posedness of the measure-valued fluid model tracking residual times.

%{\it Keywords}: { $G_t/GI/N+GI$ queues, measure-valued fluid equations, two-parameter fluid processes, elapsed service and patience times, residual service and patience times }

{\it Keywords}: {many-server queues with abandonment, time-dependent arrival rate, measure-valued fluid equations, two-parameter fluid processes, elapsed service and patience times, residual service and patience times}

\allowdisplaybreaks

\section{Introduction} \label{secIntro}

Many-server queueing models with abandonment have attracted substantial attention because of their appealing applications to customer contact centers and healthcare; see, e.g.,  \cite{B05}, \cite{GKM03}, \cite{GMR02}, \cite{G06},  
%\cite{WW02a},  
and references therein.
%  \cite{B05}, \cite{GKM03}, \cite{GMR02},   \cite{MZ04}, \cite{WW02a}, \cite{G06}, \cite{GKW07} and references therein.
% Garnett et al. (2002),  Whitt (2002), Gans et al. (2003), and Mandelbaum and Zeltyn (2004).  Brown at al. (2005).
In the $G_t/GI/N+GI$ model,  there are $N$ parallel servers, and customers arrive with a time-varying arrival rate, require i.i.d. service times,
and have i.i.d. patience times;
% abandon from the queue if their waiting times in the queue before receiving service reach their patience times, which are also assumed to be i.i.d.;
the arrival process, service and patience times are assumed to be mutually independent. The service discipline is first-come-first-served (FCFS) and non-idling, that is, no server will idle whenever there is a customer in queue.

Because of the difficulty in the exact analysis of such stochastic systems, fluid models have been recently developed to approximate the system dynamics and performance measures in a many-server heavy-traffic regime, where the arrival rate and the number of servers get large and service and patience time distributions are fixed. 
The conventional approach of using total number of customers in the system to describe system dynamics is insufficient to give a complete description and study some performance measures. Thus, measure-valued and two-parameter processes that track elapsed or residual service and patience times of each customer have been recently used  to study these stochastic models.

%Whitt \cite{WW06} pioneered to conjecture a two-parameter fluid model, which was generalized to  study the $G_t/GI/N_t+GI$ model with time-varying staffing in \cite{LW11} and \cite{LW11b}.
%Kang and Ramanan \cite{KR10} (following Kaspi and Ramanan \cite{KaR11}  for $G_t/GI/N$ models) proved a measure-valued fluid model.  Both approaches track elapsed service and patience times, but their models are developed from different perspectives. Whitt \cite{WW06} exploits the evolution equations of the density and rate functions, while Kang and Ramanan \cite{KR10} work with the dynamics of the distributions of total counts directly.  Zhang \cite{Zhang} proposed a measure-valued fluid model that tracks residual service and patience times.
%Here we show that these three fluid models are equivalent under common assumptions on the  primitives of the model.

Whitt \cite{WW06} pioneered the use of two-parameter processes to describe the system dynamics (Definition \ref{def:FM2}). In particular, $Q(t,y)$ represents the number of customers in queue at time $t$ that have waited for less than or equal to $y$, and $B(t,y)$ represents the number of customers in service at time $t$ that have received service for less than or equal to $y$.
His idea is to represent these two-parameter processes as integrals of their densities $q(t,y)$ and $b(t,y)$ with respect to $y$ (if they exist), respectively, which satisfy two fundamental evolution equations ((2.14) and (2.15) in \cite{WW06}), respectively.  
A queue boundary process  plays an important role in determining the real fluid queue size: the two-parameter density function $q(t,y)$ becomes zero for $y$ beyond the queue boundary  at each time $t$.  This approach is generalized to  study the $G_t/GI/N_t+GI$ model with both time-varying arrival rates and numbers of servers \cite{LW11} and \cite{LW11b}.

%The total service and abandonment rates at each time are written in terms of these two densities functions separately together with the associated hazard rate functions for service and patience times. A so-called queue boundary process $w(t)$ plays an important role in determining the real fluid queue size: the two-parameter density function $q(t,y)$ becomes zero for $y$ beyond the queue boundary process $w(t)$ at each time $t$. The system dynamics is then completely described by specifying the density functions $q(t,0)$ and $b(t,0)$ in the underloaded, critically loaded and overloaded regimes (equations (2.17)-(2.19) in \cite{WW06}).

%The existence and uniqueness of this two parameter fluid model follows, as a special case, from the existence and uniqueness result established in \cite{LW11, LW11b} of the two parameter fluid model for $G_t/GI/N_t+GI$ queueing model with both time-varying arrival rates and numbers of servers  under the assumptions that the system only alternates between overloaded and underloaded regimes (with a finite number of alternations in each finite time interval) and that the service and patience time distributions have piecewise continuous densities. 

Kang and Ramanan \cite{KR10}, following Kaspi and Ramanan \cite{KaR11}, used two measure-valued processes to describe the service and queueing dynamics, one
tracking the amount of time each customer has been in service, and the other tracking the amount of time each customer has spent in a potential queue, where all customers enter the potential queue upon arrival, and stay there until their patience times run out. The potential queue includes customers waiting in the real queue as well as those that have entered service or even departed but whose patience times have not run out.
% with their patience time not being reached. 
They also use a frontier waiting-time process to track the waiting time of the customer in front of the queue at each time. This frontier waiting-time process is used to determine the real fluid queue dynamics from the measure-valued process for the potential queue. The description of system dynamics is then completed by the balance equations for the fluid content processes associated with the queue, the service station and the entire system, as well as the non-idling condition; see Definition \ref{def:FMA}.

We summarize these two approaches of tracking elapsed service and patience times by stating that the two-parameter process approach in Whitt \cite{WW06} describes the system dynamics by the densities and rates, while the measure-valued process approach in Kang and Ramanan \cite{KR10} describes the system dynamics by the distributions and counting processes directly. 
The existence and uniqueness of Whitt's two-parameter fluid model are shown in discrete time under the assumption that the service and patience times have densities in \cite{WW06}. They also follow, as a special case, from the existence and uniqueness results established in \cite{LW11, LW11b} of the two-parameter fluid model for $G_t/GI/N_t+GI$ queueing model with both time-varying arrival rates and numbers of servers  under the assumptions that the system only alternates between overloaded and underloaded regimes (with a finite number of alternations in each finite time interval) and that the service and patience time distributions have piecewise continuous densities. 
The existence and uniqueness of Kang-Ramanan's fluid model are established in \cite{KR10} via the fluid limits and more recently in \cite{Kang} via the characterization of fluid model solution directly under the assumptions that the service time distribution $G^s$ has density and the hazard rate function $h^r$ of patience times is a.e. locally bounded.  Zu{\~n}iga \cite{Zu14} has recently extended Kang-Ramanan's fluid model for general service time distributions and continuous patience time distributions.

One would expect that the two approaches are equivalent since they are different formulations for the same $G_t/GI/N+GI$ queue. 
Our first main result is to establish this equivalence in Theorem \ref{wellposed}: first, a set of two-parameter fluid equations derived from the measure-valued fluid model satisfies the fluid model equations in  \cite{WW06} (see Proposition \ref{fluidtwo}), and second,  a set of measure-valued fluid equations derived from the two-parameter fluid model satisfies the fluid model equations in  \cite{KR10} (see Proposition \ref{fluidtwo2}). 
The equivalence property we establish provides a proof for the conjecture on the existence and uniqueness of Whitt's two-parameter fluid model under the assumption that the service and patience time distributions have densities (Conjecture 2.2 in \cite{WW06}). The two-parameter process formulation depends critically on the existence of the densities of the service and patience time distributions, since the densities of the two-parameter processes  may not exist for general service and patience time distributions (see Remark \ref{rmk:lem21}).

Aa a different approach, the system dynamics of $G_t/GI/N+GI$ queues can also be described by tracking residual service and patience times. It was conjectured in Section 3.3.2 of Kaspi and Ramanan \cite{KaR11} (in the case of no abandonment) that a measure-valued fluid model that tracks customers' residual service times and patience times can also be formulated in parallel to the fluid model tracking elapsed times. One advantage of considering a fluid model tracking residual times is that it enables us to easily analyze some performance measures, such as the system workload at any given time, which rely directly on the customers' residual service times; see, e.g., \cite{GW91, PW10} for infinite-server models and \cite{KaR11} for $G_t/GI/N$ queues. 
Such a fluid model tracking residual times, if suitably formulated, should be also equivalent to the above three fluid models tracking elapsed times, since all of them are formulated for the same $G_t/GI/N+GI$ queueing system. 

Zhang \cite{Zhang} provided a fluid  model tracking residual times for the $G/GI/N+GI$ model with a constant arrival rate (Definition \ref{def:FMZ}). Instead of using the potential queue as described in the fluid models tracking elapsed times, Zhang's model uses a virtual queue to describe the queueing dynamics, where all customers enter the virtual queue upon arrival and stay there until their time to enter service, which may include customers whose patience times have run out already. The existence and uniqueness of this fluid model are shown assuming continuous service time distribution and Lipschitz continuous patience time distributions \cite{Zhang}.
We study the relationship of Zhang's fluid model with the above three fluid models, in particular, focusing on Zu{\~n}iga's fluid model, and find that 
%{\color{blue}they are equivalent when the system starts from empty, but }\footnote{I would not emphasize too much on the empty initial condition. }
 they are not entirely equivalent formulations for the $G/GI/N+GI$ queue under general initial conditions; see Remarks \ref{remk:la}-\ref{remark:ari} in Section \ref{secZhanglimitation}.  
The disparity lies in the initial conditions assumed for those fluid models, in particular, the assumptions imposed on the initial contents in the virtual queue and in service in Zhang's fluid model. 
For example, in Kang-Ramanan and Zu{\~n}iga's fluid models, it is required that the residual service time of initial content in $\nu_0(dx)$ should have distribution with density $g^s(x+\cdot)/\bar G^s(x)$, whereas, in Zhang's fluid model, there is no requirement on the distribution of the residual service time of initial content in service.
 We identify the set of  necessary and sufficient conditions on the initial contents
% (see \eqref{con1} -- \eqref{DisLimRt4})
%(\ref{con1}), (\ref{con2}), (\ref{DisLimRt}) and (\ref{DisLimRt4})) 
for the equivalence of Zhang's fluid model and the above three fluid models (Theorems \ref{lem:ZKR} and \ref{lem:KRZ} and Corollary \ref{coro:KRZ}). However, in some real-life applications where the initial conditions are not satisfied (e.g., the initial content does not have the same arrival rate as new arrivals), Zhang's fluid model cannot be used.

On the other hand, from Kang-Ramanan and Zu{\~n}iga's fluid models, we obtain measure-valued fluid processes tracking residual service and patience times, which, together with the same input data as in those two fluid models, describe the service and real queueing dynamics of the same $G_t/GI/N+GI$ systems. 
%\iffalse
%The existence and uniqueness of these fluid processes tracking residual times follow from those of Zu{\~n}iga's fluid model under the general assumptions on service and patience times in \cite{Zu14}. 
%\fi 
These processes tracking residual times play an important bridging role in the discussion of the non-equivalence of  Zhang's fluid model and the fluid models tracking elapsed times.

These equivalence properties established in the paper are significant to understand the fluid dynamics of the $G_t/GI/N+GI$ model from different perspectives. 
They help to unify the different approaches in the literature, and also highlight their differences and limitations.
They provide the flexibility of choosing the most convenient approach among the different formulations, tracking elapsed or residual times, and the possibility of applying results from one formulation to another.
Some properties established with one approach can then be directly applied to other models by the equivalence relationship.
We illustrate this by two examples. 
First,  an asymptotic periodic property is proved in \cite{LW11c} for the two-parameter fluid model tracking elapsed times for the $G_t/M_t/N_t+GI_t$ queueing model, and thus, should also hold for the associated measure-valued fluid models tracking elapsed and residual times (in the special case of $G_t/M/N+GI$ queues). 
Second, it is important to show that for a fluid model, the fluid solutions converge uniformly to the steady state over all possible initial states.  That has been a difficult task for general non-Markovian many-server models. 
Thus, the equivalence property in this paper paves the way to show this with possibly any of the fluid models, whichever most convenient (see \cite{LZ} for some recent attempts in this direction).
In addition, the equivalence property results in an algorithm to compute two-parameter processes and relevant quantities under the most general conditions that cannot be computed by previous methods (see \cite{KP-alg} and its extension in \cite{MM-17} to fluid models of $G_t/GI/N+GI$ queues under the least-patient first service discipline).

Although these equivalence properties are established for the fluid limits of the associated fluid-scaled stochastic processes in the queueing model, it is conceivable that the proofs for the convergence to these fluid limits may also be unified. The two-parameter approach proves the convergence in the functional space $\sD_\sD = \sD([0,\infty), \sD([0,\infty), \RR))$ endowed with the Skorokhod $J_1$ topology. The measure-valued approach proves the convergence in the measure-valued functional space $\sD([0,\infty), \mathcal{M}([0,\infty)))$ where $\mathcal{M}([0,\infty))$ is the space of Radon measures on $\RR_+$ endowed with the Borel $\sigma$-algebra.  Tracking elapsed times enables us to use martingale arguments \cite{KR10}, but tracking residual times uses a different approach to prove the convergence \cite{Zhang}. 
%The measure-valued approach tracking residual times in \cite{Zhang} does not use any martingale convergence argument. In addition, distribution-valued processes are used to establish fluid and diffusion limits for $G/GI/\infty$ queues in \cite{DM} and \cite{RT}. These several approaches impose slightly different assumptions on the service and patience time distributions.  
So it is interesting to ask how these different approaches to establish the convergence are related and what would be the most general assumptions on the system primitives. We believe that these equivalence and coupling properties are useful in the study of other non-Markovian many-server queueing systems and networks.

\vspace{0.1in}

{\it Organization of the paper. }
The rest of the paper is organized as follows. We finish this section with some notation.  
In Section \ref{secFluidE}, we first review the definitions of the three fluid models tracking elapsed times, and then show their equivalence (Theorem \ref{wellposed}), whose proof is given in Section \ref{secConn}. 
In Section \ref{secFluidR}, we first state and discuss the fluid measure-valued processes tracking residual times derived from Kang-Ramanan and Zu{\~n}iga's fluid models in Section \ref{secMR}.
We  then review Zhang's fluid model in Section \ref{secZhang} and discuss its connection with the three fluid models tracking elapsed times in Section \ref{secZhanglimitation}. 
%Some additional proof is given in Appendix \ref{secAppProofs}. 

\vspace{0.1in}

{\it Notation. }
We use $\RR$ and $\RR_+$ to denote the spaces of real numbers and nonnegative real numbers, respectively. Given any metric space $S$, $\sC_b(S)$ is the space of bounded, continuous real-valued functions on $S$. Let $\sC_c(\RR_+)$ be the space of continuous real-valued  functions on $\RR_+$ with compact support.  
 Given a Radon measure $\xi$ on $[0,H)$ and an interval $[a,b] \subset [0,H)$, we will use $\xi[a,b]$ to denote $\xi([a,b])$. 
 Let $\mathcal D^{abs}_{[0,\infty)}(\mathcal M[0,H))$ denote the set of measure-valued processes $\mu$ with values in $\mathcal M[0,H)$, the space of Radon measures on $[0,H)$, such that for any $t\geq 0$, the measure $\int_0^t \mu_s (\cdot) ds$ is absolutely continuous with respect to the Lebesgue measure on $[0,H)$. 
Let $\mathcal D_{[0,\infty)}(\RR)$ be the space of real-valued c{\'a}dl{\'a}g functions on $[0,\infty)$. 
%The constant function $f\equiv 1$ is denoted as $\bone$.  
\iffalse
For any non-decreasing function $f$ on $\RR_+$, $f^{-1}$ denotes the inverse function of $f$ in the sense that
$
f^{-1}(y) \doteq \inf\{x\ge 0: f(x) > y\}.
$
Similarly, given a non-increasing function $f$ defined on $\RR_+$, $f^{-1}$ denotes the inverse function of $f$ in the sense that $f^{-1}(y)=\inf\{t\geq 0:\ f(t)< y\}$. 
\fi For each real-valued function $f$ defined on $[0,\infty)$, let $f^+$ and $f^-$ be the positive and the negative parts of $f$, respectively, that is, $f^+(t)=f(t)\vee 0$ and $f^-(t)=-(f(t)\wedge 0)$ for each $t\geq 0$.

\section{Fluid models tracking elapsed times} \label{secFluidE}

In the $G_t/GI/N+GI$ fluid models, we let $E(t)$ represent the cumulative amount of fluid content (representing customers) entering the system in the time interval $(0,t]$ for each $t>0$. Assume that
$E$ is a non-decreasing function defined on $[0,\infty)$ with the density function $\la(\cdot)\geq 0$, that is,
\beql{f0}
E(t)=\int_0^t \lambda(s)\,ds,\quad t\geq 0.
\eeq
%For each $t\geq 0$, $E(t)$ represents the cumulative amount of fluid content (representing customers) entering the system in the time interval $(0,t]$.
Let $G^s$ and $G^r$ denote the service and patience time distribution functions, respectively. We assume that $G^s(0+) = G^r(0+) =0$. Let
$$
H^r \doteq \inf\{x \in \RR_+: G^r(x)=1\}, \quad H^s \doteq \inf\{x \in \RR_+: G^s(x)=1\}.
$$
Then $H^r$ and $H^s$ are right supports of $G^r$ and $G^s$, respectively.

\subsection{Whitt's two-parameter fluid model} \label{secFluidtwo}

In this section we state a modified version of the two-parameter fluid model in Whitt \cite{WW06}. We assume that the functions $G^s$ and $G^r$ have density functions $g^s$ and $g^r$ on $[0,\infty)$, respectively.  Let the hazard rate functions of $G^s$ and $G^r$ be defined as $h^r \doteq g^r/\bar{G^r}$ on $[0,H^r)$ and $h^s \doteq g^s/\bar{G^s}$ on $[0,H^s)$, respectively, where $\bar{G^r} = 1- G^r$ and $\bar{G^s} = 1- G^s$. 

Let the two-parameter processes $B(t,y)$ be the amount of fluid content  in service at time $t$ that has been in service for less than or equal to $y$ units of time,  $\tilde{Q}(t,y)$ be the amount of fluid content in the potential queue at time $t$ that has been in potential queue for less than or equal to $y$ units of time, which may include the fluid content that has entered service or even departed by time $t$, and $Q(t,y)$ be the portion of $\tilde{Q}(t,y)$ that excludes the fluid content which has entered service by time $t$.
% Assume that $B(0,y) = \tilde{Q}(0,y) = Q(0,y) =0$.
Then it is obvious that $B(t, \infty)$ is the total fluid content in service and $Q(t, \infty)$ is the total fluid content in queue waiting for service.

It is assumed that these three processes are Lebesgue integrable on $[0,\infty)$ with densities $b(t,y)$, $\tilde{q}(t,y)$ and $q(t,y)$ with respect to the second component $y$, that is,
\beqal{w1}
&& B(t,y) = \int_0^y b(t,x)dx\leq 1, \quad \tilde{Q}(t,y) = \int_0^y \tilde{q}(t,x) dx\geq 0, \\
&&  Q(t,y) = \int_0^y q(t,x) dx\geq 0. \non
\eeqa
%(The existence of these densities replies on the existence of the service and patience time distribution densities.; see Remarks \ref{rmk:WW} and \ref{rmk:lem21}.)
Let $\tilde q(0,x)=q(0,x)$ as a function in $x$ have support in $[0,H^r)$ and $b(0,x)$ as a function in $x$ have support in $[0,H^s)$. 
Note that in \cite{WW06}, it is not explicitly stated that the service and patience time distributions $G^s$ and $G^r$ can be of finite support.

%the assumption of allowing finite support is not stated for the service and patience time distributions $G^s, G^r$ in \cite{WW06}.  

\begin{definition} \label{def:FM2}
A pair of functions $(B(t,y),Q(t,y))$ is a two-parameter fluid model tracking elapsed times with the input data $(\lambda(\cdot), \tilde q(0,x),b(0,x))$ if it satisfies the following conditions.

$(i)$ The service density function $b(t,x)$ satisfies 
%the first fundamental evolution equation
\beql{w2}
b(t+u, x+u) = b(t,x) \frac{\bar{G}^s(x+u)}{\bar{G}^s(x)}, \quad x \in [0,H^s), \ t \ge 0, \ u > 0.
\eeq

$(ii)$ The potential queue density function $\tilde{q}(t,x)$ satisfies 
%the second fundamental evolution equation
\beql{w3}
\tilde{q}(t+u,x+u) = \tilde{q}(t,x) \frac{\bar{G}^r(x+u)}{\bar{G}^r(x)}, \quad x \in [0,H^r), \ t\ge 0, \ u >0.
\eeq

$(iii)$ There exists a queue boundary function $w(t)$ such that  $\tilde{Q}(t,w(t))=Q(t,\infty)$ and then the queue density function $q(t,x)$ satisfies
 \beql{w5}
q(t,x) = \begin{cases}
\tilde{q}(t,x), \quad x \le w(t), \\
0, \qquad \quad x > w(t).
\end{cases}
\eeq

$(iv)$ The density functions $b(t,x)$, $\tilde{q}(t,x)$ and $q(t,x)$ satisfy the following boundary properties: 
\beql{w7}
b(t,0) = \begin{cases}
\la(t), \qifq B(t,\infty) < 1, \\
\sigma(t) \wedge \la(t), \qifq B(t,\infty) =1,  \qandq Q(t,\infty) =0,\\
\sigma(t), \qifq B(t,\infty) = 1, \qandq Q(t,\infty) >0,
\end{cases}
\eeq
\beql{w9}
\tilde q(t,0) = \lambda(t),
\eeq
and
\beql{w8}
q(t,0) = \begin{cases}
\la(t), \qifq Q(t,\infty)>0 \ (w(t)>0),\\
\la(t) - (\sigma(t) \wedge \la(t)), \qifq B(t,\infty) =1, \qandq Q(t,\infty) =0,\\
0, \qifq B(t,\infty)<1,
\end{cases}
\eeq
where \beql{w6}
\sigma(t) = \int_{[0,H^s)} b(t,x) h^s(x) dx, \quad t\ge 0.
\eeq

$(v)$ The densities $\lambda(t),\ q(t,x),\ b(t,x)$ and $\alpha(t)$ satisfy the balance equation:
\beql{w10}
\int_0^t \lambda(s) ds+\int_0^\infty q(0,x)dx= \int_0^\infty q(t,x)dx+ \int_0^t b(s,0)ds+\int_0^t\alpha(s)ds, \eeq where \beql{w6.1}
\alpha(t) = \int_{[0,H^r)} q(t,x) h^r(x) dx, \quad t\ge 0.
\eeq
\end{definition}

In \cite{WW06}, equations \eqref{w2} and \eqref{w3} are called the first and second fundamental evolution equations, respectively.
Note that the first fundamental evolution equation \eqref{w2} essentially says that the fluid content in service that has not completed service remains in service. Similarly, the second fundamental evolution equation \eqref{w3} essentially says that the fluid content in the potential queue that has not reached its patience time remains in the potential queue. For each time $t$, the queue boundary quantity $w(t)$ divides the fluid content in the potential queue into two portions. The fluid content on the left side of $w(t)$ is still in queue waiting for service and the fluid content on the right side of $w(t)$ has entered service or even departed. The quantities $b(t,0)$, $\tilde q(t,0)$, $q(t,0)$ in condition (iv) above are exactly the rates at time $t$ at which the fluid content enters service, the potential queue and the queue, respectively. The quantities $\sigma(t)$ in (\ref{w6}) and $\alpha(t)$ in (\ref{w6.1}) are precisely the total service rate and the total abandonment rate at each time $t$, respectively.
At last, the balance equation \eqref{w10} is implicit in the definition of the fluid model in \cite{WW06} and stated in equation (6) in \cite{LW11}. 

\begin{remark} \label{rmk:WW} $($Existence and uniqueness of Whitt's fluid model.$)$ Whitt \cite{WW06} has shown the existence and uniqueness of the two-parameter fluid model for $G_t/GI/N+GI$ queues in discrete time by proving a functional weak law of large numbers (FWLLN),  and conjectured them in continuous time (cf. Conjecture 2.2 of \cite{WW06}). The existence and uniqueness of the two-parameter fluid model for $G_t/GI/N_t+GI$ queues with time-dependent staffing are shown in Liu and Whitt   \cite{LW11, LW11b}, by an explicit characterization of the solution to the fluid model in \cite{LW11} and by proving an FWLLN in \cite{LW11b}, under the additional assumptions that the system only alternates between overloaded and underloaded regimes (with a finite number of alternations in each finite time interval) and that  the service and patience time distributions have piecewise continuous densities. Thus, by specializing their argument to $G_t/GI/N+GI$ queues, the conjecture is established but with the previously mentioned additional assumptions. 

In this paper, we establish the conjecture  under the assumption that the service and patience time distributions have densities, without assuming, a priori,  that the system only alternates between overloaded and underloaded regimes,  by applying the equivalence between the two fluid models in Definitions \ref{def:FM2} and \ref{def:FMA} established in Theorem \ref{wellposed} below and the existence and uniqueness of Kang-Ramanan's fluid model established in \cite{KR10, Kang}.  We remark that the existence of the densities of the service and patience time distributions is critical for the formulation of Whitt's two-parameter fluid model, because the densities of $B(t,y)$ and $Q(t,y)$ with respect to $y$ may not exist when the service and/or patience time distributions are general (see Remark \ref{rmk:lem21}).

%This is one advantage of the equivalence property: the existence and uniqueness of the two-parameter fluid model are shown without direct proofs by showing an FSLLN or FWLLN or characterizing the fluid model solution explicitly. 
\end{remark}

%The existence and uniqueness of a solution to the two-parameter fluid model are conjectured in \cite{WW06}, and shown in  \cite{LW11} and \cite{LW11b}  (with the additional feature of time-varying staffing)
%%, the existence and uniqueness of a solution to the two-parameter fluid model for $G_t/GI/N_t +GI$ queues with time-varying staffing are shown 
%by assuming that the system only alternates between overloaded and underloaded regimes (with a finite number of alternations in each finite time interval) and that  the service and patience time distributions have piecewise continuous densities. 
%We will see that the equivalence result in Theorem \ref{wellposed} below of the previously mentioned fluid models shows that the conjecture on the existence and uniqueness of a solution to the two-parameter fluid model in \cite{WW06} can be proved under the general conditions on the initial conditions, and service and patience time distributions without assuming, a priori,  that the system only alternates between overloaded and underloaded regimes. 

%\subsection{Fluid models tracking elapsed times} \label{secFluid}

\subsection{Kang-Ramanan's measure-valued fluid model} \label{secFluidmeas}

In this section, we state the measure-valued fluid model in Kang and Ramanan \cite{KR10}. They use two measure-valued processes to describe the service and queueing dynamics.  Let $\nu_t$ be a nonnegative finite measure on $[0,\infty)$ with support in $[0,H^s)$ such that $\nu_t(dx),\ x\in [0,H^s)$, represents the amount of fluid content of customers  in service whose time spent in service by time $t$ lies in the range $[x,x+dx)$.
Let $\eta_t$ be another nonnegative finite measure on $[0,\infty)$ with support in $[0,H^r)$  such that $\eta_t(dx),\ x\in [0,H^r),$ represents the amount of fluid content in the potential queue whose time spent there by time $t$ lies in the range $[x,x+dx)$, where the potential queue is an artificial queue that  includes the fluid content of customers in queue waiting for service and also the fluid content of customers that has entered service or even departed, but whose patience time has not been reached.

We assume that the functions $G^s$ and $G^r$ have density functions $g^s$ and $g^r$ on $[0,\infty)$, respectively. 
 Let $\mathcal S_0$ denote the set of triples $(\eta,\nu,x)$ such that $1 - \nu[0,H^s) = [1-x]^+$ and $\nu[0,H^s) + \eta[0,H^r)=x$, where $\eta$ is a non-negative finite measure on $[0,\infty)$ with support in $[0,H^r)$, $\nu$ is a non-negative finite measure on $[0,\infty)$ with support in $[0,H^s)$, and $x \in \RR_+$. 
The set $\mathcal{S}_0$ represents all possible measures of $(\eta,\nu)$ and values of $x$ that the initial state of the measure-valued fluid model $(\eta,\nu,X)$ can take, satisfying the non-idling condition.
\begin{definition} \label{def:FMA}
A triple of functions $(\eta,\nu,X)$ is a measure-valued fluid model tracking elapsed times with the input data $(\lambda(\cdot),\eta_0,\nu_0,X(0))$ such that $(\eta_0,\nu_0,X(0))\in \mathcal S_0$ if it satisfies the following equations. For every $\psi \in \sC_b(\RR_+)$ and $t\geq 0$,
\beql{f4}
\int_0^\infty \psi(x) \eta_t (d x) = \int_{[0,H^r)}\psi(x+t)\frac{\bar{G^r}(x+t)}{\bar{G^r}(x)}\eta_0(dx)+
\int_0^t \psi(t-s) \bar{G^r}(t-s) \la(s)ds ,
\eeq
\beql{f5}
\int_0^\infty \psi(x) \nu_t (d x) = \int_{[0,H^s)}\psi(x+t)\frac{\bar{G^s}(x+t)}{\bar{G^s}(x)}\nu_0(dx)+
 \int_{[0,t]} \psi(t-s) \bar{G^s}(t-s) d K(s),
\eeq
where \begin{eqnarray}
K(t) &=& B(t) + D(t) - B(0) = \nu_t [0,H^s)+ D(t)- \nu_0 [0,H^s),  \label{33} \\ 
D(t) &=& \int_0^t \left( \int_{[0,H^s)} h^s(x) \nu_s(dx) \right) ds , \label{f80}
\end{eqnarray}
\begin{eqnarray}
E(t) +Q(0)&=& Q(t) + K(t) + R(t),  \label{31}
\end{eqnarray}
\beql{f8}
R(t) = \int_0^t \left( \int_{[0,\chi(s)]} h^r(x) \eta_s(dx) \right) ds,
\eeq
 \beql{f9}
\chi(s) = \inf\{x\in [0,H^r): \eta_s[0,x] \ge Q(s)\},
\eeq
\beql{f1}
Q(t) = (X(t) - 1)^+, 
\eeq
\beql{f2}
B(t) =  \nu_t [0,\infty) = X(t) \wedge 1 = 1 - (1-X(t))^+, 
\eeq
and
\beql{f12}
Q(t)(1-B(t)) = 0.
\eeq
\end{definition}

%It is shown in \cite{KR10} that the above fluid model has a unique solution  under the assumptions that the hazard rate functions $h^r$ and $h^s$  are either bounded or lower semi-continuous. Recently a new direct proof of the existence and uniqueness of the fluid model is established in \cite{Kang} under the assumptions that the service time distribution $G^s$ has density and the hazard rate function $h^r$ is a.e. locally bounded.  

\iffalse
\begin{remark}
When the system starts empty and the arrival rate $\lambda(\cdot)$ is constant, the assumption on the existence of densities of $G^s$ and $G^r$  can be relaxed in the formulation of the fluid model since the processes $D$ and $R$ have the following alternative representations:
\[R(t)=\int_0^t G^r(\chi(s))ds.\] and
\[D(t) = \int_0^t \left( \int_0^s K(s-u) dG^s(u) \right) ds.\]
\end{remark}
\fi

In this fluid model, $B(t)$ represents the total fluid content of customers in service, $Q(t)$ represents the total fluid content of customers in queue waiting for service, and $X(t)$ represents the total fluid content of customers in the system at each time $t$. Then, by (\ref{f1}) and (\ref{f2}),  \beql{32}
X(t)=B(t)+Q(t).
\eeq  The additional quantities $K(t),\ R(t),\ D(t),\ \chi(t)$ can naturally be interpreted, respectively, as the cumulative amount of fluid content that has entered service by time $t$, the cumulative amount of fluid content that has abandoned from the queue by time $t$, the amount of fluid content that has departed the system after service completion by time $t$, and the waiting time of the fluid content at the head of the queue at time $t$, that is, the fluid content in queue with the longest waiting time.

 For completeness, we now provide an intuitive explanation for these fluid equations.  The equation (\ref{f4}) governs the evolution of the measure-valued process $\eta_t$.  Note that when $x\leq t$, the amount of fluid content $\eta_t(dx)$ is the fraction of the amount of fluid content $\lambda(t-x)$ arriving to the system at time $t-x$ and whose time in the system since its arrival is more than $x$ by time $t$. It is easy to see that this fraction equals to $\bar{G^r}(x)$. When $x>t$, the amount of fluid content $\eta_t(dx)$ is the fraction of the amount of fluid content $\eta_0(d(x-t))$ initially in queue and whose waiting time is more than $x$ by time $t$ given that it is more than $x-t$ at time $0$. This fraction equals to $\bar{G^r}(x)/\bar{G^r}(x-t)$. This shows that (\ref{f4}) holds. A similar observation yields (\ref{f5}). The equations (\ref{33})--(\ref{31})  are simply mass conservation equations for the queue and the server station, respectively. Since $\nu_s(dx),\ x\in [0,s],$ represents the amount of fluid content in service whose time in service lies in the range $[x,x+dx)$ at time s, and $h^s(x)$ represents the fraction of the amount of fluid content with time in service $x$ (that is, with service time no less than $x$) that would depart from the system while having time in service in $[x, x + dx)$. Hence, it is natural to expect $\int_{[0,H^s)}h^s(x)\nu_s(dx)$ to represent the departure rate of fluid content from the fluid system at time $s$ and thus, expect (\ref{f80}) holds. A similar explanation can be applied to (\ref{f8}) except that, to consider the real reneging rate, we can only consider $x<\chi(s)$ since all the fluid content with the time in the system more than $\chi(s)$ has entered service by time $s$. The equation (\ref{f12}) represents the usual non-idling condition.
 
 By adding (\ref{31}) and (\ref{33}) together and using (\ref{32}), we see that
\beql{total}
E(t) + X(0)=X(t)+R(t)+D(t).
\eeq
By the representations of $E$, $R$ and $D$ in (\ref{f0}), (\ref{f8}) and (\ref{f80}), we have from (\ref{total}) that $X$ is absolutely continuous. In turn, using the fact that $|[n- a]^+ - [n -b]^+| \leq |a - b|$, it is easy to see from (\ref{f2}) and (\ref{33}) that $B$ and then $K$ are absolutely continuous. So there exists a Lebesgue integrable function $\kappa$ such that
\beql{Kcont}
K(t)=\int_0^t \kappa(s)ds,\quad t\ge 0.
\eeq
By (\ref{33}) and (\ref{f80}), the process $K$ has the following representation:
 \beql{t20}
K(t) = B(t)-B(0) + \int_0^t \left( \int_{[0,H^s)}  h^s(x) \nu_s(dx)\right) ds.
\eeq
Then it follows from  the same argument as in deriving (3.12) of \cite{KaR11} that the process $\kappa$ satisfies for a.e. $t\in \RR_+$,
\beql{kappa}
\kappa(t)=\left\{\begin{array}{ll} \lambda(t) & \mbox{ if } X(t)<1, \\
\lambda(t)\wedge \int_{[0,H^s)} h^s(x)\nu_t(dx) & \mbox{ if } X(t)=1, \\ \int_{[0,H^s)} h^s(x)\nu_t(dx)  & \mbox{ if } X(t)>1. \end{array}\right.
\eeq

\begin{remark} \label{rem:euo} $($Existence and uniqueness of Kang-Ramanan's fluid model.$)$
Under the assumptions that the hazard rate functions $h^r$ and $h^s$  are either bounded or lower semi-continuous, Kang and Ramanan \cite{KR10} established the existence of the measure-valued fluid model in Definition \ref{def:FMA} by proving an FWLLN and also showed its uniqueness via  the fluid model characterization. The existence and uniqueness of Kang-Ramanan's fluid model directly from the characterization of its solution is established in Kang \cite{Kang}, under the weaker assumptions that the service time distribution $G^s$ has density and the hazard rate function $h^r$ is a.e. locally bounded. 
\end{remark}

Now we state our first result on the equivalence between the two fluid models  described in Definitions \ref{def:FM2} and \ref{def:FMA}.   Its proof is deferred to Section  \ref{secConn}. As a consequence, it also gives a proof for Conjecture 2.2 of \cite{WW06} under the assumption that the service and patience time distributions have densities and $h^r$ is a.e. locally bounded.
%we have the following theorem which establishes Conjecture 2.2 of \cite{WW06}.
\begin{theorem} \label{wellposed}
Existence and uniqueness of Whitt's fluid model in Definition \ref{def:FM2} is equivalent to existence and uniqueness of Kang-Ramanan's fluid model in Definition \ref{def:FMA} for the $G_t/GI/N+GI$ queue with the time-dependent arrival rate $\la(\cdot)$ and the initial data $(\eta_0,\nu_0,X(0))\in \mathcal S_0$, where $\eta_0(dx)=\tilde q(0,x)dx=q(0,x) dx$ and $\nu_0(dx)=b(0,x)dx$.
\end{theorem}

\subsection{Zu{\~n}iga's fluid model} \label{secFluidmeasE}
 Recently, Zu{\~n}iga \cite{Zu14} extended Kang-Ramanan's fluid model without assuming that the patience time distribution $G^r$ and service time distribution $G^s$ have densities. In this section, we state this extended Kang-Ramanan's fluid model and establish some useful properties on certain quantities in the model, which are needed in the subsequent analysis.

Define a measure $M^r$ on $[0,H^r]$ by \[dM^r(x)\doteq \bone_{\{x<H^r\}}\bar{G}^r(x-)^{-1} d G^r(x) + \bone_{\{G^r(H^r-)<1\}}\delta_{H^r}(dx),\] and a measure $M^s$ on $[0,H^s]$ by \[dM^s(x)\doteq \bone_{\{x<H^s\}}\bar{G}^s(x-)^{-1} d G^s(x) + \bone_{\{G^s(H^s-)<1\}}\delta_{H^s}(dx).\]

\begin{definition} \label{def:EFMA}
% $($Zu{\~n}iga's fluid model in  \cite{Zu14}$)$
A triple of processes $(\eta,\nu,X)\in \mathcal D^{abs}_{[0,\infty)}(\mathcal M[0,H^r))\times \mathcal D^{abs}_{[0,\infty)}(\mathcal M[0,H^s))\times \mathcal D_{[0,\infty)}(\RR)$ is a solution to an extended Kang-Ramanan's measure-valued fluid model with the input data $(\lambda(\cdot),\eta_0,\nu_0,X(0))$ such that $(\eta_0,\nu_0,X(0))\in \mathcal S_0$ if $q_t$ and $p_t$, the densities of $\int_0^t \nu_s (\cdot) ds$ and $\int_0^t \eta_s (\cdot) ds$, respectively,  satisfy the following conditions. 
There exist $K(\cdot)$, a process of bounded variation started at $0$, $\chi(\cdot)$, $B(\cdot)$, $Q(\cdot)$, $D(\cdot)$, $R(\cdot)$ such that for every $\psi \in \sC_b(\RR_+)$ and $t\geq 0$, \eqref{f4}--\eqref{33}, \eqref{31}, \eqref{f9}--\eqref{f12} hold and  \begin{eqnarray}
D(t) &=& \int_{[0,H^s]} q_t(x) dM^s(x) , \label{Ef80}
\end{eqnarray}
\beql{Ef8}
R(t) = \int_{[0,H^r]} \int_{[0,t]} \bone_{\{x\leq \chi(s)\}} d_s p_s(x) dM^r(x),
\eeq
where the integral with respect to $p_s(x)$ is defined  as a Lebesgue-Stieltjes integral in $s$. 
\end{definition}

\begin{remark} \label{rem:new}
Zu{\~n}iga's fluid model stated in Definition \ref{def:EFMA} is equivalent to Definition 3.4 of \cite{Zu14} due to Lemma 4.1 and Remark 4.2 of \cite{Zu14} and the given input data $(\lambda(\cdot),\eta_0,\nu_0,X(0))$. The main difference of Zu{\~n}iga's fluid model from Kang-Ramanan's fluid model in Definition \ref{def:FMA} is that the processes $D$ and $R$ satisfy \eqref{Ef80} and \eqref{Ef8} instead of \eqref{f80} and \eqref{f8} due to the lack of existence of densities of $G^s$ and $G^r$, respectively. 
By Lemma 4.1 of \cite{Zu14}, the densities $q_t$ and $p_t$ can be written as
\beql{qt-Zu}
q_t(x) = \bar{G}^s(x-) K((t-x)^+) + \int_{[(x-t)^+, x)} \frac{\bar{G}^s(x-)}{\bar{G}^s(y)} \nu_0(dy),
\eeq
and
\beql{pt-Zu}
p_t(x) = \bar{G}^r(x-) E((t-x)^+) + \int_{[(x-t)^+, x)} \frac{\bar{G}^r(x-)}{\bar{G}^r(y)} \eta_0(dy).
\eeq
When $G^s$ and $G^r$ are assumed to have densities, $g^r$ and $g^s$, respectively, Zu{\~n}iga's fluid model is reduced to Kang-Ramanan's fluid model. Zu{\~n}iga's fluid model admits a unique solution (established in Theorem 3.5  via an FWLLN and Theorem 4.4 via the characterization of the fluid model in \cite{Zu14}) under the assumptions that $G^r$ is continuous, $\eta_0$ is diffuse, and $\nu_0$ is diffuse if $G^s$ is not continuous (Assumption 3.1 of \cite{Zu14}). 
\end{remark}

We end this section by showing the following {\color{red} critical }lemma for Zu{\~n}iga's fluid model in Definition \ref{def:EFMA}, which will be used in Section \ref{secFluidR} in discussing the relationship of the fluid models tracking elapsed times stated in this section and a fluid  model tracking residual times stated in Section \ref{secZhang}.
%Its proof  is given in Appendix \ref{secAppProofs}. 

\begin{lemma}\label{lem:chiinc}
In Definition \ref{def:EFMA}, the processes  $D$ and $R$ have the following representations: for each $t\geq 0$,
%\beql{mr621}
%B(t) 
%=\int_{[0,H^s)} \frac{\bar G^s(y+t)}{\bar G^s(y)} \nu_0(dy) +
%\int_0^t  \bar G^s(t-s) dK(s),
%\eeq
\beql{mrd}
D(t) = \int_{[0,H^s)} \frac{G^s(y+t)-G^s(y)}{\bar G^s(y)} \nu_0(dy) +
 \int_0^t G^s(t-s) dK(s) ,
\eeq
%\beql{mr61}
%Q(t) = 1_{\{\chi(t)\ge t \}}\int_{[0,\chi(t)- t ]} \frac{\bar G^r(y+t)}{\bar G^r(y)}\eta_0(dy) + 
%\int_{(t-\chi(t))^+}^t \bar G^r(t-s) \la(s)ds ,
%\eeq
\beqal{mrq}
R(t) &=& 
\int_{[0,H^r)} \left(\int_{(y,y+t]} \bone_{\{y\leq \chi(x-y)-(x-y)\}}  dG^r(x)\right) \bar{G}^r(y)^{-1}\eta_0(dy)  \non\\
&&  \qquad + \int_0^t \int_{[0,H^r]} \bone_{\{x\leq s \wedge \chi(s)\}} \la(s-x) dG^r(x) ds . 
\eeqa
Moreover, the process $K(t)$ is non-decreasing and the process $\chi(t)$ satisfies the following property:
 \begin{eqnarray} \label{chiinc} \chi(t)-\chi(s)\leq t-s \mbox{ whenever } 0\leq s<t<\infty. \end{eqnarray} 
\end{lemma}

\begin{remark} \label{rmk:lem21}
It is evident that the representation of the process $D$ in \eqref{mrd} implies that $D(t)$ is not absolute continuous when the service time distribution does not have density. 
Thus, we cannot write the total service rate (departure rate) as in \eqref{w6}. 
Although the two-parameter processes $B(t,y)$, $\tilde{Q}(t,y)$ and $Q(t,y)$ can be obtained as in \eqref{k2} from the Zu{\~n}iga's fluid model $(\nu_t, \eta_t, X)$ in Definition \ref{def:EFMA}, their densities with respect to $y$ may not exist and the associated two-parameter fluid model using densities $b(t,x)$ and $q(t,x)$ cannot be formulated with the densities as in Definition \ref{def:FM2}. 
\end{remark}

%\section{Proofs of Lemma 2.1.}
%\label{secAppProofs}
\noindent {\it Proof of Lemma \ref{lem:chiinc}.}
%{\it Proof.} 
By \eqref{Ef80} and \eqref{qt-Zu}, applying interchange of the order of integration and integration by parts, we easily obtain \eqref{mrd}. 
%\beqal{mrdp1}
%D(t) 
%%&=&  \int_{[0, H^s]} q_t(x) d M^s(x) \non\\
% &=& \int_{[0, H^s]} \Bigg(  \int_{x-t}^x \frac{\bar{G}^s(x-)}{\bar{G}^s(y)} \nu_0(dy) + \bar{G}^s(x-) K(t-x)  \Bigg) d M^s(x) , \non \\
%&=& \int_{[0,H^s]} \frac{G^s(y+t)-G^s(y)}{\bar G^s(y)} \nu_0(dy) +
% \int_0^t G^s(t-s) dK^\ell(s) \non
%\eeqa
To show $R(t)$ in \eqref{mrq}, from (\ref{Ef8}) and (\ref{pt-Zu}), we obtain that for each $t\geq 0$,
\begin{eqnarray} \label{Rr1}
 && \ R(t) = \int_{[0,H^r]} \int_0^t \bone_{\{x\leq \chi(s)\wedge s\}} \bar{G}^r(x-)\la(s-x)ds dM^r(x)   \\ && \quad  +  \int_{[0,H^r]} \int_{[0,t]} \bone_{\{x\leq \chi(s)\}} d_s \left(\int_{[(x-s)^+,x)} \frac{\bar{G}^r(x-)}{\bar{G}^r(y)} \eta_0(dy)\right) dM^r(x) \non \\ &=&  \int_0^t \int_{[0,H^r]} \bone_{\{x\leq \chi(s)\wedge s\}} \la(s-x)\bar{G}^r(x-)dM^r(x) ds \non \\ && \ + \int_{[0,H^r]} \int_{[0,x\wedge t]} 1_{\{x\leq \chi(s)\}} d_s \left(\int_{[x-s,H^r)} \bone_{\{y<x\}}\frac{\bar{G}^r(x-)}{\bar{G}^r(y)} \eta_0(dy)\right) dM^r(x) \non \\ &=& \int_0^t \int_{[0,H^r]} \bone_{\{x\leq \chi(s)\wedge s\}} \la(s-x)\bar{G}^r(x-)dM^r(x) ds  \non \\ && \ + \int_{[0,H^r]} \int_{[[x-t]^+,x]} \bone_{\{x\leq \chi(x-s)\}} \bone_{\{s<x\}}\frac{\bar{G}^r(x-)}{\bar{G}^r(s)} \eta_0(ds) dM^r(x), \non \\
&=& \int_0^t \int_{[0,H^r]} \bone_{\{x\leq s \wedge \chi(s)\}} \la(s-x)\bar{G}^r(x-) dM^r(x) ds \non \\ & & \ + \int_{[0,H^r)} \left(\int_{(y,y+t]} \bone_{\{y\leq \chi(x-y)-(x-y)\}} \bar{G}^r(x-) dM^r(x)\right) \bar{G}^r(y)^{-1}\eta_0(dy) \non \\ &=& \int_0^t \int_{[0,H^r]} \bone_{\{x\leq s \wedge \chi(s)\}} \la(s-x) dG^r(x) ds \non \\ & & \ + \int_{[0,H^r)} \left(\int_{(y,y+t]} \bone_{\{y\leq \chi(x-y)-(x-y)\}} dG^r(x)\right) \bar{G}^r(y)^{-1}\eta_0(dy) , \non
\end{eqnarray}
where the second term in the second equality  follows from Theorem 3.6.1 of \cite{Bo} with $X=[[x-t]^+,x]$ $Y=[0,x\wedge t]$, $f(s)=x-s$ and $\mu$ such that $\mu[a,b] = \int_{[a,b)}\frac{\bone_{\{y<x\}}}{\bar{G}^r(y)} \eta_0(dy)$ and the last equality follows from the interchange of the order of integrations.

We next prove the non-decreasing property of $K(t)$. 
It follows from this representation of $R(t)$ in \eqref{mrq} that Lemma 4.4 of \cite{KR10} holds, that is, for any $0 \le a \le b < \infty$, if $Q(t) =0$ (equivalently, $\chi(t) =0$) for all $t \in [a,b]$, then $R(b) - R(a) =0$. Then the proof for the non-decreasing property of $K(t)$ will follow the same argument in Lemma 4.5 in \cite{KR10} using \eqref{Rr1}.

We now prove the property of $\chi(t)$ in \eqref{chiinc}. 
By a similar argument as in Lemma 3.4 of \cite{KR10} on time shifts, to prove the lemma, without loss of generality, we may assume that $s=0$ in (\ref{chiinc}). Suppose that the property of $\chi(t)$ in \eqref{chiinc} does not hold, that is, there is a time $t_2> 0$ such that $\chi(t_2)>\chi(0)+t_2$. Let \[t_1 \doteq \sup\{u\leq t_2:\ \chi(u)\leq \chi(0)+u\}.\] 
Then $\chi(t_1-)\leq \chi(0)+t_1$ and for each $u\in [t_1,t_2]$, \beql{chiu} \chi(u)\geq  \chi(0)+u\geq \chi(t_1-)+(u-t_1) \mbox{ and } \chi(t_2)>\chi(t_1-)+(t_2-t_1). \eeq 
By \eqref{Rr1}, it is clear that $R(t)-R(t-)\geq 0$ for each $t>0$. By applying the above display and time shift at $t_1$, we have 
\begin{eqnarray*}
& & R(t_2)-R(t_1) \\ &= & \int_{[0,H^r)} \left(\int_{(y,y+t_2-t_1]} \bone_{\{y\leq \chi(t_1+x-y) -(x-y)\}} dG^r(x)\right) \bar{G}^r(y)^{-1}\eta_{t_1}(dy) \\ & & \qquad +
 \int_0^{t_2-t_1} \left(\int_{[0,H^r]}\bone_{\{u\leq s\wedge \chi(t_1+s)\}}  \la(t_1+s-u)dG^r(u) \right)ds.
\end{eqnarray*}
 It follows from (\ref{chiu}) that $s\wedge \chi(t_1+s)=s$ and $\chi(t_1+s)-s\geq  \chi(t_1-)$ for each $s\in (0,t_2-t_1]$. Hence the above display implies that 
 \begin{eqnarray*}
& & R(t_2)-R(t_1-) \\ &\geq &
\int_{[0,H^r)} \left(\int_{(y,y+t_2-t_1]} \bone_{\{y\leq \chi(t_1-)\}} dG^r(x)\right) \bar{G}^r(y)^{-1}\eta_{t_1}(dy) \\ & & \qquad +
 \int_0^{t_2-t_1} \left(\int_{[0,H^r]}\bone_{\{u\leq s\}}  \la(t_1+s-u)dG^r(u) \right)ds \\ &=& \int_{[0,H^r)} \bone_{\{y\leq \chi(t_1-)\}} \frac{G^r(y+t_2-t_1) - G^r(y)}{\bar{G}^r(y)}\eta_{t_1}(dy) \\ & & \quad + \int_{0}^{t_2-t_1} G^r(t_2-t_1-u) \la(t_1+u)du, \non  
\end{eqnarray*}
where 
and the second term on the right hand side of the last display  follows from Proposition 0.4.5 of \cite{RM}.  
Since (\ref{31}) holds for Zu{\~n}iga's fluid model $(\eta,\nu,X)$ and $K$ is non-decreasing, then  the above three displays imply that  
\begin{eqnarray*}
Q(t_2) &= & Q(t_1-) + (E(t_2)-E(t_1-)) -(R(t_2)-R(t_1-)) - (K(t_2)-K(t_1-)) \\ &\le & \eta_{t_1}[0,\chi(t_1-)] +\int_{t_1}^{t_2} \la(u)du - \int_{0}^{t_2-t_1} G^r((t_2-t_1-u)-) \la(t_1+u)du \\ & & \quad -\int_{[0,H^r)}\bone_{[0,\chi(t_1-)]}(x)\frac{G^r(x+(t_2-t_1))-G^r(x)}{\bar{G^r}(x)}\eta_{t_1}(dx) \\ &=& \int_{0}^{t_2-t_1} \bar{G^r}(t_2-t_1-u) \la(t_1+u)du \\ & & \quad + \int_{[0,H^r)}\bone_{[0,\chi(t_1-)+(t_2-t_1)]}(x+(t_2-t_1))\frac{\bar{G^r}(x+(t_2-t_1))}{\bar{G^r}(x)}\eta_{t_1}(dx) \\ &=& \eta_{t_2}[0,\chi(t_1-)+(t_2-t_1)]. 
\end{eqnarray*}
From this and the definition of $\chi$, we have $\chi(t_2)\leq \chi(t_1-)+(t_2-t_1)$, which contradicts (\ref{chiu}). Thus, the lemma is proved. 
~~~\bsq

\section{Measure-valued fluid models tracking residual times} \label{secFluidR}

%Tracking customers' elpased service times and patience times, as described in Kang-Ramanan's fluid model, is not the only approach to capture system dynamics of a fluid many-server queueing system. It was conjectured in Section 3.3.2 of Kaspi and Ramanan \cite{KaR11} (in the case of no abandonment) that a measure-valued fluid model that tracks customers' residual service times and patience times can also be derived in parallel to the fluid model tracking elapsed times. One advantage of considering a fluid model tracking residual times is that it enables us to easily analyze the performance measures, such as the system workload at any given time, that rely directly on the customers' residual service times. Such a fluid model tracking residual times, if it exists, should be closely connected with the fluid model tracking elapsed times in Definition \ref{def:FMA}, since both fluid  models are formulated for the same  many-server queueing system. 

%Zhang \cite{Zhang} provided a measure-valued fluid model for the $G/GI/N+GI$ model by tracking residual service and patience times of each customer when the arrival rate is constant and the initial condition for the fluid model satisfies certain conditions.  

We first state the two measure-valued processes tracking residual times that arise from Zu{\~n}iga's fluid model for the same $G_t/GI/N+GI$ queueing system in Section \ref{secMR}. We then state Zhang's fluid model in Section \ref{secZhang}, and discuss its connection with the three fluid models tracking elapsed times in Section \ref{secZhanglimitation}. The two measure-valued processes tracking residual times introduced in Section \ref{secMR} play an important bridging role in making the connection. 
%In particular, we emphasize the advantages and disadvantages of Zhang's fluid model. 
%To overcome the disadvantages, we present a new fluid model tracking residual times in Section \ref{secnewmodel}, which is shown to be equivalent to the two fluid models tracking elapsed times in Section \ref{secnewequiv}.

\subsection{Measure-valued processes tracking residual times from Zu{\~n}iga's fluid model}
\label{secMR} 
\iffalse
Let $\nu^\ell_t$  be a nonnegative finite measure on $[0,\infty)$ with support in $[0,H^s)$ such that $\nu^\ell_t(dx)$ represents the amount of fluid content of customers in service whose residual service time at time $t$ lies in the range $[x,x+dx)$. Similarly, let $\eta^\ell_t$  be a nonnegative finite measure on $[0,\infty)$ with support in $[0,H^r)$ such that $ \eta^\ell_t(dx)$ represents the amount of fluid content of customers in queue whose residual patience time at time $t$ lies in the range $[x,x+dx)$. 
\fi
Zu{\~n}iga's fluid model naturally give rise to the following two measure-valued processes $\nu^\ell_t$ and 
 $\eta^\ell_t$. 
% that keep track of residual times of fluid contents of customers in queue and in service.
For each $t\geq 0$, clearly the mapping
\begin{eqnarray*}\psi &\mapsto & \int_{[0,H^s)}\left(\int_{(y+t,\infty)} \frac{\psi(x-y-t)}{\bar G^s(y)}d G^s(x)\right)\nu_0(dy) \\
&& \qquad + \int_{[0,t]} \left(\int_{(t-s,\infty)} \psi(x-t+s)d G^s(x)\right)dK(s)\end{eqnarray*} is a positive linear functional on $\mathcal{C}_c(\RR_+)$ since $K$ is non-decreasing by Lemma \ref{lem:chiinc}. Then by Riesz-Markov-Kakutani representation theorem, there is a unique regular Borel measure $\nu^\ell_t$ with  support $[0,H^s)$ such that for every $\psi \in \mathcal{C}_b(\RR_+)$, 
\beqal{mr2}
\int_0^\infty \psi(x) \nu^\ell_t(dx) &=&  \int_{[0,H^s)}\left(\int_{(y+t,\infty)} \frac{\psi(x-y-t)}{\bar G^s(y)}d G^s(x)\right)\nu_0(dy) \non\\
&& \qquad + \int_0^{t} \left(\int_{(t-s,\infty)} \psi(x-t+s)d G^s(x)\right)dK(s).
\eeqa
Similarly, for each $t\geq 0$, there is a unique regular Borel measure $\eta^\ell_t$ with support $[0,H^r)$ such that for every $\psi \in \mathcal{C}_b(\RR_+)$, 
\beqal{mr4}
\int_0^\infty \psi(x)\eta^\ell_t(dx) & = &  \bone_{\{\varsigma(t)\leq 0\}}\int_{[0,-\varsigma(t)]}\left(\int_{(y+t,\infty)} \frac{\psi(x-y-t)}{\bar G^r(y)}dG^r(x)\right)\eta_0(dy) \non\\ && \qquad + \int_{\varsigma^+(t)}^{t} \left(\int_{(t-s,\infty)} \psi(x-t+s) d G^r(x) \right)\lambda(s) ds ,
\eeqa
where \beql{chiw}
\varsigma(t) = t - \chi(t).  
\eeq
Since $\chi(t)$ represents the elapsed patience time of the fluid content of customers that has been in queue the longest at time $t$, then the quantity $\varsigma(t)$ can be interpreted as the arrival time of the fluid content of customers that has been in queue the longest at time $t$.  It is clear that $ \varsigma(t)\leq t$ for each $t\ge 0$. 
At time 0, $\varsigma(0) = - \chi(0)$ represents the arrival time of the oldest fluid content in queue initially, and thus, it follow from \eqref{f9} that
\beql{mrvs1}
\varsigma(0)=-\inf\{x\in [0,H^r):\ \eta_0[0,x) \geq X(0)-\nu_0[0,H^s)\} .
\eeq
 
%If $\varsigma(t) <0$ at time $t> 0$, it means that the initial fluid content of customers at time 0 has not all entered service yet by time $t$, and can be interpreted as the arrival time of the initial fluid content that has been in queue the longest at time $t$. 

We first argue that $\nu^\ell$ and $\eta^\ell$ are two measure-valued processes tracking residual times of fluid content of customers in service and in queue, respectively. 
\iffalse
 So it is clear that
\beql{mrvs2}
 \varsigma(t)\leq t  \quad \text{for each} \  t \ge 0,
\eeq
 and intuitively, we see that $\varsigma(t)$ can be related to $\chi(t)$  as in \eqref{chiw}.
At time 0, $\varsigma(0) = - \chi(0)$ represents the arrival time of the oldest fluid content in queue initially, and thus, it follow from \eqref{f9} that
\beql{mrvs1}
\varsigma(0)=-\inf\{x\in [0,H^r):\ \eta_0[0,x) \geq X(0)-\nu_0[0,H^s)\} .
\eeq
 If $\varsigma(t) <0$ at time $t> 0$, it means that the initial fluid content of customers at time 0 has not all entered service yet by time $t$, and can be interpreted as the arrival time of the initial fluid content that has been in queue the longest at time $t$. 
\fi

For each $z\geq 0$, by plugging $\psi(x)=\bone_{(z,\infty)}(x)$ into \eqref{mr2} and \eqref{mr4}, we have 
\begin{eqnarray}
\nu^\ell_t(z,\infty) &=& \int_{[0,H^s)}\frac{\bar{G}^s(y+t+z)}{\bar G^s(y)}\nu_0(dy)  \label{nuz} \\
&& \qquad \qquad + \int_{[0,t]} \bar{G}^s(t-s+z)dK(s),\non \\ 
\eta^\ell_t(z,\infty) &=& \bone_{\{\varsigma(t)\leq 0\}}\int_{[0,-\varsigma(t)]}\frac{\bar{G}^r(y+t+z)}{\bar G^r(y)}\eta_0(dy) \label{etaz}\\
&& \qquad \qquad+ \int_{\varsigma^+(t)}^{t} \bar{G}^r(t-s+z)\lambda(s) ds. \non
\end{eqnarray}
By \eqref{f4} and \eqref{f5}, we obtain
\begin{eqnarray*}
\eta_{t+z}[z,\chi(t)+z] &=& \int_{[0,H^r)}\bone_{[z,\chi(t)+z]}(y+t+z)\frac{\bar{G^r}(y+t+z)}{\bar{G^r}(y)}\nu_0(dy) \\ & & \quad +
 \int_0^t\bone_{[z,\chi(t)+z]}(t+z-s)\bar{G^r}(t+z-s) \la(s) ds \\ &=& \bone_{\{\varsigma(t)\leq 0\}}\int_{[0,-\varsigma(t)]}\frac{\bar{G}^r(y+t+z)}{\bar G^r(y)}\eta_0(dy) \\
&& \quad + \int_{\varsigma^+(t)}^{t} \bar{G}^r(t-s+z)\lambda(s) ds ,
\end{eqnarray*}
and
\begin{eqnarray*}
\nu_{t+z}[z,\infty) =\int_{[0,H^s)}\frac{\bar{G^s}(y+t+z)}{\bar{G^s}(y)}\nu_0(dy)+
 \int_{[0,t]}\bar{G^s}(t+z-s) d K(s).
\end{eqnarray*}
Hence,  we obtained the following {\it coupling} property between $(\nu, \eta)$ and $(\nu^\ell, \eta^\ell)$:
\beql{mrme1}
\nu^\ell_t(z,\infty) = \nu_{t+z}[z,\infty) \mbox{ and } \eta^\ell_t(z,\infty) = \eta_{t+z}[z,\chi(t)+z],\ quad z\geq 0. 
\eeq
Intuitively, $\nu_{t+z}[z,\infty)$ represents the amount of fluid content in service at time $t+z$ with elapsed service time at least $z$, which is precisely the amount of fluid content in service at time $t$ that will still be in service at time $t+z$  and then is equal to the amount of fluid content in service at time $t$ that has residual service time greater than $z$. (Note that the fluid content in service at time $t$ that has residual service time exactly equal to $z$ will depart from service and hence will not be in service at time $t+z$.) Thus, by the first equality in  (\ref{mrme1}), $\nu^\ell_t(z,\infty)$  represents the amount of fluid content in service at time $t$ that has residual service time greater than $z$, that is, $\nu^\ell_t$ keeps track of the residual time of fluid content in service at time $t$. Similarly, $\eta_{t+z}[z,\chi(t)+z]$ represents the amount of fluid content in the potential queue at time $t+z$ with elapsed patience time between $z$ and $\chi(t)+z$, which is precisely the amount of fluid content in queue at time $t$ that will not abandon by time $t+z$. This amount of fluid content is equal to the amount of fluid content that has residual patience time more than $z$ units of time at time $t$ and then is represented by $\eta^\ell_t(z,\infty)$ by the second equality in (\ref{mrme1}). Then $\eta^\ell_t$ keeps track of the residual patience times of customers in queue at time $t$.

When $t=0$, \eqref{nuz}, \eqref{etaz} and  \eqref{mrme1}  become: for each $z\geq 0$,
\beqal{mr22}
\nu^\ell_0(z,\infty) = \int_{[0,H^s)}\frac{\bar{G}^s(y+z)}{\bar G^s(y)}\nu_0(dy) = \nu_z [z,\infty),
\eeqa
\beql{mr42}
\eta^\ell_0(z,\infty) = \bone_{\{\varsigma(0)\leq 0\}}\int_{[0,-\varsigma(0)]}\frac{\bar{G}^r(y+z)}{\bar G^r(y)}\eta_0(dy) = \eta_z[z,\chi(0)+z].
\eeq

\begin{remark} 
%The fluid model in Definition \ref{def:FMR} does not require  the patience time distribution $G^r$ and service time distribution $G^s$ to have densities, respectively. 
When $G^r$ and  $G^s$ have densities $g^r$ and $g^s$, respectively,  \eqref{mr2} and \eqref{mr4} are equivalent to the following representations:
\beqal{mr21}
\int_0^\infty \psi(x) \nu^\ell_t(dx) &=&  \int_{[0,H^s)}\left(\int_0^\infty \frac{g^s(y+t+x)}{\bar G^s(y)}\psi(x)dx\right)\nu_0(dy) \non\\
&& \qquad + \int_{[0,t]} \left(\int_0^\infty g^s(t-s+x)\psi(x)dx\right)dK(s),
\eeqa
\beqal{mr41}
\int_0^\infty \psi(x)\eta^\ell_t(dx) & = &  \bone_{\{\varsigma(t)\leq 0\}}\int_{[0,-\varsigma(t)]}\left(\int_0^\infty \frac{g^r(y+t+x)}{\bar G^r(y)}\psi(x)dx\right)\eta_0(dy) \non\\ && \qquad + \int_{\varsigma^+(t)}^{t} \left(\int_0^\infty g^r(t-s+x) \psi(x)dx \right)\lambda(s) ds.
\eeqa
In this case, for each $t\geq 0$, the two measures $\eta^\ell_t$ and $\nu^\ell_t$ have
densities $b_\ell(t,x)$ and $q_t(t,x)$, respectively, which can be expressed as
\beql{tr56}
b_\ell(t,y) = \int_{[0,H^s)}\frac{g^s(x+t+y)}{\bar G^s(x)}\nu_0(dx) + \int_{[0,t]} g^s(y+t-u) dK(u),
\eeq
and 
\beql{tr567}  q_\ell(t,y) = \bone_{\{\varsigma(t)\leq 0\}}\int_{[0,-\varsigma(t)]}\frac{g^r(x+t+y)}{\bar G^r(x)}\eta_0(dx)+ \int_{\varsigma^+(t)}^{t} g^r(y+t-u) \lambda (u) du.
\eeq
\iffalse
In addition to $(\nu^\ell,\eta^\ell)$,  we also introduce a measure-valued process $\tilde{\eta}^\ell$ satisfying that for every $\psi \in \mathcal{C}_b(\RR_+)$, 
\beqal{mr62}
\int_0^\infty \psi(x) \tilde{\eta}^\ell_t(dx) &=&  \int_{[0,H^r)}\left(\int_0^\infty \frac{g^r(y+t+x)}{\bar G^r(y)}\psi(x)d x\right)\eta_0(dy) \non\\
&& \qquad + \int_0^{t} \left(\int_0^\infty g^r(t-u+x)\psi(x)dx\right)\lambda(u) du.
\eeqa 
For each $t\geq 0$ and $x\geq 0$, $\tilde{\eta}^\ell_t(dx)$ represents the amount of fluid content of customers in the potential queue whose residual patience times at time $t$ lie in the range $[x,x+dx)$.
For each $t\geq 0$, the measure $\tilde \eta^\ell_t$ admits a density function $\tilde{q}_\ell(t,y)$ with the representation:
\beql{mr64}
 \tilde{q}_\ell(t,y)=\int_{[0,H^r)}\frac{g^r(x+t+y)}{\bar G^r(x)}\eta_0(dx)+ \int_0^{t} g^r(y+t-s) \lambda(s) ds.
\eeq
%Moreover, by applying interchange of the order of integration,  (\ref{mrq}) is equivalent to the following.
%\beqal{mrq1}
%R^\ell(t) &=& \int_0^t \left(\int_0^{\varsigma^-(s)}\frac{g^r(y+s)}{\bar G^r(y)}\eta_0(dy)\right)ds \\ & & \qquad +
% \int_0^t \left(\int_{\varsigma^+(s)}^s  g^r(s-u) \la(u)d u \right)ds. \non
%\eeqa
\fi
\end{remark}

\subsection{Zhang's  fluid model} \label{secZhang}

 Zhang \cite{Zhang} uses a so-called {\it virtual queue} to describe the queueing dynamics, instead of the potential queue used in the three fluid models in Section \ref{secFluidE}. 
In the definitions of both potential and virtual queues, all customers enter them upon arrival. 
{\it  The difference between them lies in how customers depart.}
 Customers can leave the potential queue only when their patience expires, that is, at the instant when their remaining patience times are zeros. 
Whereas, customers can only leave the virtual queue in their turns of service.
Customers in the virtual queue may have already run out of patience (i.e., the remaining patience time is negative) at their turns of service. Whenever a server becomes free, the server will check the oldest customer in the virtual queue. If the customer being checked has not abandoned yet (its remaining patience time is still positive), then the server will start serving this customer and this customer is removed from the virtual queue, and otherwise, this customer is simply removed from the virtual queue and the server will turn to check the next oldest customer. 
We now state Zhang's fluid model.

\begin{definition} \label{def:FMZ}
$($Zhang's fluid model in \cite{Zhang}.$)$
Assume that the fluid arrival rate $\lambda(t) = \la $ for each $t\ge 0$, where $\lambda> 0$ is a constant. A pair of measure-valued processes $(\mathcal{R}, \mathcal{Z})$ is a solution to the fluid model if the following conditions are satisfied:

$(i)$ $(\mathcal{R}, \mathcal{Z})$ satisfies the following two equations:
\beql{r1}
\mathcal{R}_t(C_x) = \la \int_{t- Q_v(t)/\la}^t \bar{G}^r(t+x-s) ds, \quad x \in \RR,
\eeq
and
\beql{r2}
\mathcal{Z}_t(C_x) = \mathcal{Z}_0(C_x+t) + \int_0^t \bar{G}^r(Q_v(s)/\la) \bar{G}^s(t+x-s) d L_v(s), \quad x \in \RR_+,
\eeq where $C_x \doteq (x,\infty)$ for $x \in \RR$, $Q_v(t) = \mathcal{R}_t(\RR)$ is of bounded variation  and $L_v(t) = \lambda t - Q_v(t)$;

$(ii)$ the non-idling conditions in \eqref{f2} and \eqref{f12} hold for $B(t)=\mathcal{Z}_t(\RR_+)$, $Q(t)=\mathcal{R}_t(\RR_+)$ and $X(t)=B(t)+Q(t)$; 

$(iii)$ the initial condition $(\mathcal{R}_0,\mathcal{Z}_0)$ satisfies
\beql{r3}
\mathcal{R}_0(C_x) = \la \int_0^{ Q_v(0)/\la}\bar{G}^r(x+s) ds, \ x \in \RR, \mbox{ and }\mathcal{Z}_0(\{0\}) = 0, 
\eeq
and the non-idling condition at time $0$ in \eqref{f2} and \eqref{f12}.

\end{definition}

\vspace{0.1in}
In Zhang's fluid model, $\mathcal{R}_t(C_x)$ can be interpreted  as the fluid content of customers in the virtual queue with residual patience times strictly bigger than $x$ and $\mathcal{Z}_t(C_x)$ can be interpreted as the fluid content of customers in service with residual service times  strictly bigger than $x$ at each time $t$.
Then $Q(t)$, $B(t)$, $Q_v(t)$ and $L_v(t)$ represent, respectively, the total fluid content of the real queue at time $t$,  the total fluid content of customers in service at time $t$,  the total fluid content in the virtual queue at time $t$, and the cumulative customers removed from the virtual queue by time $t$. 
The existence and uniqueness of Zhang's fluid model are proved in Theorem 3.1 of \cite{Zhang} by an explicit characterization of its solution, 
 under the assumptions that the service time distribution $G^s$ is continuous and the patience time distribution $G^r$ is Lipschitz continuous.  
%It is worth pointing out that in Zhang's fluid model, the patience time distribution $G^r$ and the service time distribution $G^s$ are not required to have densities, respectively. This is mainly due to the arrival rate being assumed constant and the special structure of  $\mathcal{R}_0$  in (\ref{r3}).

\subsection{Connection between Zhang's fluid model and the three fluid models in Section \ref{secFluidE}} \label{secZhanglimitation}
Among the three fluid models in Section \ref{secFluidE}, we have showed in Theorem \ref{wellposed} that Whitt's fluid model is equivalent to Kang-Ramanan fluid model  and in Remark \ref{rem:new} that Zu{\~n}iga's fluid model extends Kang-Ramanan's fluid model by relaxing the assumption on the existence of densities of $G^r$ and $G^s$. Since Zhang's fluid model keeps track of customers' residual times and does not need $G^r$ and $G^s$ to have densities (\cite{Zhang} does assume that $G^s$ is continuous and $G^r$ is Lipschitz continuous  to establish existence and uniqueness), while Zu{\~n}iga's fluid model keeps track of customers' elapsed times and also does not need $G^r$ and $G^s$ to have densities, it is natural to question if Zhang's fluid model and Zu{\~n}iga's fluid model are in fact  equivalent in describing system dynamics of the same $G_t/GI/N+GI$ queues. If so, this will enable researchers to borrow results from either one of the two to study the system performance of $G_t/GI/N+GI$ queues. 

In this section we provide a detailed discussion on Zhang's fluid model in connection with Zu{\~n}iga's fluid model (and hence Kang-Ramanan's fluid model and Whitt's fluid model). 
The three fluid models in Section \ref{secFluidE} allow time-varying arrival rate $\lambda(\cdot)$, whereas, Zhang's fluid model requires a constant arrival rate $\lambda$. Thus the discussion in this section will focus on the four formulations with a constant arrival rate. 
We first show by a series of remarks that {\it Zhang's fluid model is not entirely equivalent to the three fluid models tracking elapsed times for the same $G/GI/N+GI$ queueing system under general initial conditions, that is, Zhang's fluid model and the three fluid models tracking elapsed times may not be formulated simultaneously for the same $G/GI/N+GI$ queueing system under certain general initial conditions.}  

\begin{remark} $($On the arrival rate.$)$ \label{remk:la}
The imposed condition on $\mathcal{R}_0$ in Zhang's fluid model requires that the initial fluid content of customers in the virtual queue depends on the arrival rate $\la$ after time $0$, whereas in real life applications, the customers' arrival patterns before time $0$ and after time $0$ are likely different. 
%hence there is no reason to impose this requirement.  
Thus, Zhang's fluid model may not be appropriate for those applications.
 %with $G/GI/N+GI$ queueing systems. 
 In contrast, the three fluid models tracking elapsed times do not have this restriction. 
\end{remark}

\begin{remark}\label{remk:z} $($The initial condition on $\mathcal{R}_0$.$)$
Zhang's fluid model requires that the system initial condition $\mathcal{R}_0$ satisfies \eqref{r3}, that is, 
\beql{zR0x}
\mathcal{R}_0(C_x) = \la \int_0^{ Q_v(0)/\la}\bar{G}^r(x+s) ds, \ x \in \RR.
\eeq
Let $\mathcal{R}_0^+$ be the restriction  of $\mathcal{R}_0$ on $[0,\infty)$. 
Then $\mathcal{R}_0^+$  keeps track of the residual patience times of the fluid content of customers initially in queue.  So if Zhang's fluid model were equivalent to Zu{\~n}iga's fluid model for the same  $G/GI/N+GI$ queueing system assuming a constant arrival rate, we must have 
$\mathcal{R}_0^+ = \eta^\ell_0$ in \eqref{mr42}, that is,
%It will be shown in Section  \ref{secnewmodel}  (see \eqref{??} when $t=0$), that 
$$
\mathcal{R}_0^+(C_x)=\bone_{\{\varsigma(0)\leq 0\}}\int_{[0,-\varsigma(0)]}\frac{\bar{G}^r(y+x)}{\bar G^r(y)}\eta_0(dy), \quad  \forall x\geq 0,
$$
%$$
%  \int_{[0,H^r)}\frac{\bar{G}^r(y+x)}{\bar G^r(y)}\eta_0(dy),\quad  \forall x\geq 0,
%$$
where $\eta_0$ is the initial condition for the $\eta$ in Definition \ref{def:EFMA}, and then 

\beql{DisLimR0}
 \la \int_0^{ Q_v(0)/\la}\bar{G}^r(x+s) ds =\bone_{\{\varsigma(0)\leq 0\}}\int_{[0,-\varsigma(0)]}\frac{\bar{G}^r(y+x)}{\bar G^r(y)}\eta_0(dy),\quad  \forall x\geq 0.
\eeq
%\beql{DisLimR0}
% \la \int_0^{ Q_v(0)/\la}\bar{G}^r(x+s) ds =  \int_{[0,H^r)}\frac{\bar{G}^r(y+x)}{\bar G^r(y)}\eta_0(dy),\quad  \forall x\geq 0.
%\eeq
We first note that there may not be a unique $\eta_0$ satisfying \eqref{DisLimR0} for the given $\mathcal{R}_0$. For example, when $G^r$ has density $g^r(x)= e^{-x},\ x\in \RR_+$, 
$$
\bone_{\{\varsigma(0)\leq 0\}}\int_{[0,-\varsigma(0)]}\frac{\bar{G}^r(y+x)}{\bar G^r(y)}\eta_0(dy) 
	= \bone_{\{\varsigma(0)\leq 0\}} e^{-x}\eta_0[0,-\varsigma(0)], 
$$
%$$\int_{[0,H^r)}\frac{\bar{G}^r(y+x)}{\bar G^r(y)}\eta_0(dy)= e^{-x} \eta_0[0,H^r)$$
 and  
  $$\la \int_0^{ Q_v(0)/\la}\bar{G}^r(x+s) ds = \la e^{-x} \left(1-e^{-Q_v(0)/\la} \right).
  $$ Thus, any $\eta_0$ satisfying $ \bone_{\{\varsigma(0)\leq 0\}}\eta_0[0,-\varsigma(0)]= \la (1-e^{-Q_v(0)/\la})$  will satisfy \eqref{DisLimR0}.

Moreover, it is clear that the above display \eqref{DisLimR0} does not hold for an arbitrary initial condition $\eta_0$. 
For example, if $\eta_0(dx)= \lambda^\dag \bar G^r(x)dx$ for some positive $\lambda^\dag \neq \lambda$, then 
\begin{eqnarray*}
&& \bone_{\{\varsigma(0)\leq 0\}}\int_{[0,-\varsigma(0)]}\frac{\bar{G}^r(y+x)}{\bar G^r(y)}\eta_0(dy) 
	\\
&=& \lambda^\dag \bone_{\{\varsigma(0)\leq 0\}}\int_{[0,-\varsigma(0)]} \bar{G}^r(y+x)dy ,
% \neq \mathcal{R}_0^+(C_x) ,
\end{eqnarray*}
which is not equal to $\mathcal{R}_0^+(C_x)$ in \eqref{zR0x} even if $-\varsigma(0) = Q_v(0)/\la$. 
%\[ \int_{[0,H^r)}\frac{\bar{G}^r(y+x)}{\bar G^r(y)}\eta_0(dy)= \lambda^\dag \int_{[0,H^r)}\bar{G}^r(y+x)dy \neq \mathcal{R}_0^+(C_x).\]   
 Thus, for a fluid $G/GI/N+GI$ queueing system with a constant arrival rate $\lambda$ after time $0$, the initial conditions $\eta_0(dx)= \lambda^\dag \bar G^r(x)dx$ for $\lambda^\dag \neq \lambda$ and $(\nu_0,X(0))$ such that $(\eta_0,\nu_0,X(0))\in \mathcal S_0$,  Zu{\~n}iga's fluid model can be well formulated, but there is no corresponding Zhang's fluid model $(\mathcal{R}, \mathcal{Z})$ that describes the same system. 
\end{remark}

\begin{remark} \label{remark:ari} $($The initial condition on $\mathcal{Z}_0$.$)$
Zhang's fluid model only requires that $\mathcal{Z}_0(\{0\}) = 0$. This condition is rather general. We show by an example that for a $G/GI/N+GI$ queueing system, although Zhang's fluid model can be formulated with that initial condition $\mathcal{Z}_0$,  there may not exist an (unique) initial measure $\nu_0$ to formulate a corresponding Zu{\~n}iga's fluid model for the same system. 

Consider the service time distribution $G^s$ being exponential with unit rate, that is, $g^s(x)= e^{-x},\ x\in \RR_+$. Let $\mathcal{Z}_0$ be the measure that tracks the residual service times of fluid content of customers initially in service and satisfies $\mathcal{Z}_0(\{0\}) = 0$, and assume that Zhang's fluid model can be formulated with $\mathcal{Z}_0$.  Suppose that  Zu{\~n}iga's fluid model can also be formulated for some measure $\nu_0$, which tracks the elapsed service times of fluid content of customers initially in service.  
%It is shown in Section  \ref{secnewmodel} (see equation ???) that 
%It was suggested in Section 3.3.2 of Kaspi and Ramanan \cite{KaR11}, and 
By \eqref{mr22}, 
 if Zhang's fluid model and Zu{\~n}iga's fluid model were equivalent, $\mathcal{Z}_0$ and $\nu_0$ must satisfy the following equation:
\[ \mathcal{Z}_0(C_x) = \int_{[0,H^s)}\frac{\bar G^s(y+x)}{\bar G^s(y)}\nu_0(dy)= \nu_0[0,H^s) e^{-x}, \quad x\geq 0. \] 
If the given $\mathcal{Z}_0$ satisfies $ \mathcal{Z}_0(C_x) = c e^{-x}$ for some constant $c>0$, then any such measure $\nu_0$ in Zu{\~n}iga's fluid model satisfying $c =\nu_0[0,H^s)$ will satisfy the above display. However, on the other hand,  if the given $\mathcal{Z}_0$, satisfying $\mathcal{Z}_0(\{0\}) = 0$, does not have an exponential density, then this contradicts the above equation resulting from the equivalence, and implies that no corresponding measure $\nu_0$ can be found for Zu{\~n}iga's fluid model to be well formulated for the given queueing system. 

\end{remark}

From the discussion in Remarks \ref{remk:la}, \ref{remk:z} and \ref{remark:ari}, it is clear that the class of fluid many-server queueing systems where Zhang's fluid model can be formulated is \emph{not} the same as the class of fluid many-server queueing systems where Zu{\~n}iga's fluid model (and hence Kang-Ramanan's fluid model and Whitt's fluid model) can be formulated. 

We next look more closely into the conditions on fluid $G/GI/N+GI$ queueing systems where  Zhang's fluid model and the three fluid models tracking elapsed times can all be used to describe the  system dynamics for the same system. To simplify the exposition, we focus on Zhang's fluid model and Zu{\~n}iga's fluid model.  
Our findings are stated in the following two theorems.  

\begin{theorem} \label{lem:ZKR}
Given a Zhang's fluid model $(\mathcal{R}, \mathcal{Z})$ for a $G/GI/N+GI$ queueing system with arrival rate $\la$, there exists a Zu{\~n}iga's fluid model $(\eta,\nu,X)$ for the same queueing system with the input data $(\la,\eta_0,\nu_0,X(0))$ such that $(\eta_0,\nu_0,X(0))\in \mathcal S_0$ with 
 \begin{equation} \label{con4}
\eta_0(dx) \doteq  \la \bone_{[0,Q_v(0)/\la]}(x) \bar{G}^r(x)dx ,
\end{equation}
if and only if,  for the given $\mathcal{Z}_0$, $\nu_0$ satisfying \begin{eqnarray} 
\mathcal{Z}_0(C_x) = \int_{[0,H^s)}\frac{\bar{G}^s(y+x)}{\bar G^s(y)}\nu_0(dy), \ x\geq 0. \label{con2}\end{eqnarray}
% In particular, we can choose 
% \begin{equation} \label{con4}
%\eta_0(dx) \doteq  \la \bone_{[0,Q_v(0)/\la]}(x) \bar{G}^r(x)dx .
%\end{equation}
\end{theorem}
\noindent {\it Proof.} The ``only if" part follows directly from the discussion in Remark \ref{remark:ari}. We now focus on ``if" part. 

Let $\eta_0$ be as given in (\ref{con4}). For each $t\geq 0$, the following mapping
\begin{eqnarray*}\psi &\mapsto & \int_{[0,H^r)}\psi(x+t)\frac{\bar{G^r}(x+t)}{\bar{G^r}(x)}\eta_0(dx)+
\la \int_0^t \psi(t-s) \bar{G^r}(t-s) ds\end{eqnarray*} is a positive linear functional on $\mathcal{C}_c(\RR_+)$. Then by Riesz-Markov-Kakutani representation theorem, there is a unique regular Borel measure $\eta_t$ on $\RR_+$ such that (\ref{f4}) holds. It is clear that $\eta_t$ has support $[0,H^r)$. 

For each $t\geq 0$, define 
% $\eta_0= \la 1_{[0,Q_v(0)/\la]}(x) \bar{G}^r(x)dx$,
 \[K(t)\doteq \int_0^t \bar{G}^r(Q_v(s)/\la) d L_v(s)\mbox{ and }  R(t)\doteq \lambda 
 \int_0^t G^r(Q_v(s)/\la)\,ds. \]
Then, for each $t\geq 0$, with the above $K$ and the given $\nu_0$ satisfying (\ref{con2}), the mapping
\begin{eqnarray*}\psi &\mapsto & \int_{[0,H^s)}\psi(x+t)\frac{\bar{G^s}(x+t)}{\bar{G^s}(x)}\nu_0(dx)+
 \int_{[0,t]} \psi(t-s) \bar{G^s}(t-s) d K(s) \end{eqnarray*} is a positive linear functional on $\mathcal{C}_c(\RR_+)$. By Riesz-Markov-Kakutani representation theorem, there is a unique regular Borel measure $\nu_t$ that satisfies (\ref{f5}). Let $B,\ Q,\ X$ be the associated processes in Zhang's fluid model and for each $t\geq 0$, define \[D(t)\doteq \int_{[0,H^s)} \frac{G^s(y+t)-G^s(y)}{\bar G^s(y)} \nu_0(dy) +
 \int_0^t G^s(t-s) dK(s).\] We show that $(\eta,\nu,X)$ satisfies Definition \ref{def:EFMA}. 

From (\ref{r1}), it is clear that 
\beql{qqv}
 Q(t)= \mathcal{R}_t(\RR_+) = \la \int_{t-Q_v(t)/\la}^t \bar{G}^r(t-s) ds =\la G^r_d(Q_v(t)/\lambda),\eeq where $G^r_d(x)=\int_0^x \bar G^r(s)ds$.
It is established in the proof of Theorem 3.1 of \cite{Zhang} that $Q(t)/\lambda < G^r_d(\infty)= G^r_d(H^r)$. Then it follows that  $Q_v(t)/\la <H^r$. Since $Q_v$ is of bounded variation by (\ref{r1}),  it follows that $Q$ is also of bounded variation and by the chain rule formula (Proposition 4.6 in Chapter 0 of \cite{RM}) \[Q(t) = Q(0)+ \int_0^t \bar G^r(Q_v(s)/\lambda)dQ_v(s).\]
Thus, by the definition of $K$ and the above display for $Q(t)$,\begin{eqnarray*} K(t) &=&  Q(0)-\left(Q(t)- \la\int_0^t \bar G^r(Q_v(s)/\lambda) ds\right) \\ &=& Q(0)-\left(Q(t)- \la\int_0^t \bar G^r((G^r_d)^{-1}(Q(s)/\lambda)) ds\right). \end{eqnarray*}
Then it follows from Lemma A.3 of \cite{Zhang} that $K$ is non-decreasing. 
Simple calculation also shows that \begin{eqnarray*} & & Q(t)+K(t)+R(t) \\ &=& Q(0)+ \int_0^t \bar G^r(Q_v(s)/\lambda)dQ_v(s) + \int_0^t \bar{G}^r(Q_v(s)/\lambda) d (\lambda s -Q_v(s)) \\ & & \qquad  + \lambda 
 \int_0^t G^r(Q_v(s)/\la)\,ds \\ &=& Q(0)+\lambda t, \end{eqnarray*}
which establishes (\ref{31}). 
 For each $t\geq 0$, define $\chi(t)$ by the right hand side of (\ref{f9}). It follows from the construction of $\eta_t$ and the given $\eta_0$ in \eqref{con4} that 
\begin{eqnarray} \label{QQv1}
 Q(t) &=& \eta_t[0,\chi(t)] = \la \int_0^{[\chi(t)-t]^+\wedge Q_v(0)/\la}\bar{G^r}(x+t)dx \non\\
&& \qquad \qquad \qquad + \la
\int_{[t-\chi(t)]^+}^t \bar{G^r}(t-s) ds.
 \end{eqnarray}
When $\chi(t)>t$, the above display is reduced to 
$$
Q(t)/\la  = \int_0^{\chi(t)\wedge (t+Q_v(0)/\la)} \bar{G}^r(s)ds. 
$$
 Comparing this with (\ref{qqv}), we have  $Q_v(t)/\la = \chi(t)\wedge (t+Q_v(0)/\la).
$
%\beql{Qvchi}
%Q_v(t)/\la = \chi(t)\wedge (t+Q_v(0)/\la).
%\eeq 
When $\chi(t)\leq t$, the  display in \eqref{QQv1} is reduced to $Q(t)/\la  = \int_0^{\chi(t)} \bar{G}^r(s)ds$ and hence $Q_v(t)/\la = \chi(t)$. Combining the two cases, we have for each $t\geq 0$, 
\beql{Qvchi1}
 Q_v(t)/\la = \chi(t)\wedge (t+Q_v(0)/\la). 
\eeq
%From this, it is easy to check that $\eta_0$ defined in (\ref{con4}) satisfies (\ref{con1}), (\ref{DisLimRt}) and \eqref{DisLimRt4}.
% {\color{red} need to also assume that $\chi(0) = Q_v(0)/\la$??}

For each $t\geq 0$, it follows from \eqref{con4} and the definition of $M^r$ that
\begin{eqnarray*}
&& \int_{[0, H^r]} \int_{[0,t]} \bone_{\{x \le \chi(s) \}} d_s \left( \int_{[(x-s)^+, x)} \frac{\bar{G}^r(x-)}{\bar{G}^r(y)}   \eta_0(dy )  \right)  d M^r(x)  \\ & & \quad +\la \int_{[0, H^r]} \int_{[0,t]} \bone_{\{x \le \chi(s) \wedge s\}} \bar{G}^r(x-) ds  d M^r(x) \\
&=&  \la \int_{[0, H^r]} \int_{[0,x\wedge t]} \bone_{\{x \le \chi(s)\}} d_s \left( \int_{[x-s, x)} \frac{\bar{G}^r(x-)}{\bar{G}^r(y)}    \bone_{[0,Q_v(0)/\la]}(y) \bar{G}^r(y)dy  \right)  d M^r(x) \\ & & \quad + \la  \int_{[0,t]} \int_{[0, H^r]} \bone_{\{x \le \chi(s) \wedge s\}}  d G^r(x) ds  \\
&=&  \la \int_{[0, H^r]} \int_0^{t\wedge x} \bone_{\{x \le \chi(s)\}} d_s \left( \int_{[x-s, x)}    \bone_{[0,Q_v(0)/\la]}(y) dy  \right)  d G^r(x) \\ & & \quad + \la  \int_0^t G^r(\chi(s) \wedge s) ds  \\
&=&  \la \int_{[0, H^r]} \int_0^{t} \bone_{\{s\leq x \le \chi(s) \}}  \bone_{[0,Q_v(0)/\la]}(x-s) ds    d G^r(x) \\ & & \quad + \la  \int_0^t G^r(\chi(s) \wedge s) ds  \\
&=&  \la  \int_0^t  G^r(\chi(s)\wedge (s+Q_v(0)/\la)) ds =   \la  \int_0^t  G^r(Q_v(s)/\la) ds ,
\end{eqnarray*}
where the last equality follows from \eqref{Qvchi1}. This,  together with the definition of $R(t)$ and  \eqref{pt-Zu}, implies that \eqref{Ef8} holds. 

By using (4.5) of \cite{Zu14}, we obtain
\begin{eqnarray*}
 &&\int_{[0, H^s]} \left( \bar{G}^s(x-) K([t-x]^+) + \int_{[x-t]^+}^x \frac{\bar{G}^s(x-)}{\bar{G}^s(y)} \nu_0(dy) \right) d M^s(x) \\
&=& \int_{[0,H^s)} \frac{G^s(y+t)-G^s(y)}{\bar G^s(y)} \nu_0(dy) +
 \int_0^t G^s(t-s) dK(s) 
\end{eqnarray*}
which is equal to the process $D(t)$ by definition, and implies that \eqref{Ef80} holds. 

For each $t\geq 0$, (\ref{33}) holds by applying interchange of the order of integration to (\ref{f5}) and using the definitions of $D$ and $B$.  The properties  (\ref{f1})--(\ref{f12}) follow from property (ii) of Zhang's fluid model. Thus, this completes the proof that $(\eta,\nu,X)$ is a Zu{\~n}iga's fluid model satisfying Definition \ref{def:EFMA}. Clearly from the construction, both the given Zhang's fluid model and the constructed Zu{\~n}iga's fluid model  describe the same $G/GI/N+GI$ queueing system. ~~~\bsq

\begin{theorem} \label{lem:KRZ}
Given  a Zu{\~n}iga's fluid model $(\eta,\nu,X)$  for a $G/GI/N+GI$ queueing system with the input data $(\la,\eta_0,\nu_0,X(0))$ such that $(\eta_0,\nu_0,X(0))\in \mathcal S_0$, there exists a Zhang's fluid model $(\mathcal{R}, \mathcal{Z})$ for the same queueing system
with arrival rate $\la$ if and only if  $\eta_0$ satisfies the following condition: for each $t\geq 0$, there exists a solution $z_t$, independent of $x\geq 0$, to the equation in $z$:
\beql{DisLimRtCond}  
 \la \int^{z}_{t\wedge \chi(t)} \bar{G}^r(x+s) ds 
=  \bone_{\{\chi(t)\geq t\}}\int_{[0,\chi(t)-t]}\frac{\bar{G}^r(y+t+x)}{\bar G^r(y)}\eta_0(dy),
\eeq
such that \begin{eqnarray} 
 \la \int_0^t G^r(z_s)ds  
&=&   \la \int_0^t G^r(\chi(s)\wedge s) ds \non\\ & & \quad +  \int_{[0,H^r)} \left(\int_{(y,y+t]} \bone_{\{x\leq \chi(x-y)\}} dG^r(x)\right) \bar{G}^r(y)^{-1}\eta_0(dy). \label{DisLimRtCond1} 
\end{eqnarray}
In this case, $\mathcal{Z}_0$ can be chosen as defined by \eqref{con2} for the given $\nu_0$, and $\mathcal{R}_0$ can be chosen as defined by \eqref{r3} for $Q_v(0)=z_0 \la$, where $z_0$ is the solution, independent of $x\geq 0$, that satisfies \eqref{DisLimRtCond}  for $t=0$. 
\end{theorem}

\begin{remark}
When $G^r$ has a density $g^r$, the conditions \eqref{DisLimRtCond} and \eqref{DisLimRtCond1} can be replaced as follows: for each $t\geq 0$, there exists a solution $z_t$, independent of $x\geq 0$, to the equation in $z$: 
\beql{DisLimRtCond2}  
\la G^r(x+z) 
= \la G^r(x+t\wedge \chi(t))+ \bone_{\{\chi(t)\geq t\}}\int_{[0,\chi(t)-t]}\frac{g^r(y+t+x)}{\bar G^r(y)}\eta_0(dy).
\eeq
 In fact, \eqref{DisLimRtCond} follows from \eqref{DisLimRtCond2} directly by integrating both sides of \eqref{DisLimRtCond} in $x$. 
It follows from \eqref{DisLimRtCond2} with $x=0$ that \begin{eqnarray*}
& &  \la \int_0^t G^r(\chi(s)\wedge s) ds \non\\ & & \quad +  \int_{[0,H^r)} \left(\int_{(y,y+t]} \bone_{\{x\leq \chi(x-y)\}} dG^r(x)\right) \bar{G}^r(y)^{-1}\eta_0(dy) \\ &=&  \la \int_0^t G^r(\chi(s)\wedge s) ds \non\\ & & \quad +  \int_{[0,H^r)} \left(\int_0^t \bone_{\{y\leq \chi(x)-x\}} g^r(x+y)dx\right) \bar{G}^r(y)^{-1}\eta_0(dy) \\ &=& \la \int_0^t G^r(\chi(s)\wedge s) ds + \int_0^t \la(\bar{G}^r(s\wedge \chi(s))-\bar{G}^r(z_s)) ds \\ &=& \la \int_0^t G^r(z_s)ds.
\end{eqnarray*}
Thus, \eqref{DisLimRtCond1} holds.
\end{remark}

\noindent {\it Proof of Theorem \ref{lem:KRZ}.}  We first show the ``only if" part. Recall that in Zhang's fluid model,  $\mathcal{R}_t^+$, the restriction  of $\mathcal{R}_t$ on $[0,\infty)$, tracks the residual patience times of the fluid content of customers in queue at time $t$. 
If there exists a Zhang's fluid model $(\mathcal{R}, \mathcal{Z})$ to describe the same  $G/GI/N+GI$ queueing system together with Zu{\~n}iga's fluid model $(\eta,\nu,X)$,  the measure  $\mathcal{R}_t$ must satisfies  (see \eqref{etaz})  
$$ \mathcal{R}_t(C_x)=\bone_{\{\chi(t)\geq t\}}\int_{[0,\chi(t)-t]}\frac{\bar{G}^r(y+t+x)}{\bar G^r(y)}\eta_0(dy) + \la\int_{0}^{t\wedge \chi(t)}  \bar{G}^r(s+x) ds, $$ 
for each $t\ge 0$ and $x\ge0$, and 
% So if Zhang's fluid model $(\mathcal{R}, \mathcal{Z})$ and Kang-Ramanan's fluid model $\eta,\nu,X)$ can be used to describe the same  $G/GI/N+GI$ queueing system with the constant arrival rate $\la$, 
hence $\eta_0$ must satisfy that for each $t\ge 0$ and $x\geq 0$, 
\beql{DisLimRt}
 \la \int^{Q_v(t)/\la}_{t\wedge \chi(t)} \bar{G}^r(x+s) ds 
=  \bone_{\{\chi(t)\geq t\}}\int_{[0,\chi(t)-t]}\frac{\bar{G}^r(y+t+x)}{\bar G^r(y)}\eta_0(dy). 
\eeq
When $t=0$, (\ref{DisLimRt}) is reduced to 
\beql{con1}
 \mathcal{R}_0(C_x) =  \bone_{\{\chi(0)\geq 0\}}\int_{[0,\chi(0)]}\frac{\bar{G}^r(y+x)}{\bar G^r(y)}\eta_0(dy),\  x\geq 0,
\eeq
which is discussed in Remark \ref{remk:z}.
 Moreover, in Zhang's fluid model, since customers in queue will renege when their residual patience times reach zero, then by differentiating (\ref{r1}) in $x$ and letting $x=0$, we have the abandonment rate at time $t$ is given by \[\la \left(\bar{G}^r(x) - \bar{G}^r(x + Q_v(t)/\la)\right)\mid_{x=0}= G^r(Q_v(t)/\la).
\]
Then $R(t)$, the cumulative abandonment by time $t$, is given by $\int_0^t G^r(Q_v(s)/\la)ds$. On the other hand, by Zu{\~n}iga's fluid model, $R(t)$ is given by (\ref{mrq}). Then $\eta_0$ must also satisfy that for each $t>0$, 
\begin{eqnarray}
&& \quad  \la \int_0^t G^r(Q_v(s)/\la)ds  =  \la \int_0^t G^r(\chi(s)\wedge s) ds   \label{DisLimRt4}  \\ & & \qquad \qquad \qquad +  \int_{[0,H^r)} \left(\int_{(y,y+t]} \bone_{\{x\leq \chi(x-y)\}} dG^r(x)\right) \bar{G}^r(y)^{-1}\eta_0(dy).  \non 
\end{eqnarray}
Note that for each $t\geq 0$, $Q_v(t)$ satisfies \eqref{DisLimRtCond} and \eqref{DisLimRtCond1}, independent of $x\geq 0$. Hence the ``only if" part is established.

For the ``if" part, let $\mathcal{Z}_0$ and $\mathcal{R}_0$ be defined as in the statement of the theorem. It is clear that the defined $\mathcal{R}_0$ and $\eta_0$ satisfy (\ref{con1}) and $(\mathcal{R}_0, \mathcal{Z}_0)$ satisfies property (iii) of Zhang's fluid model.
Let $\chi(t)$, $B(t)$, $Q(t)$, $K(t)$, $D(t)$, $R(t)$ be the associated auxiliary processes from Zu{\~n}iga's fluid model $(\eta,\nu,X)$.
For each $t>0$, define \[Q_v(t)\doteq \la z_t\ \mbox{ and } L_v(t) \doteq  \lambda t - Q_v(t),\] where $z_t$ is the solution, independent of $x\geq 0$, that satisfies (\ref{DisLimRtCond}) and (\ref{DisLimRtCond1}).
%\[Q_v(t)\doteq \inf\left\{y\geq 0:\ \la \int_0^y \bar{G}^r(s)ds \geq Q(t)\right\}\la \mbox{ and } L_v(t) \doteq  \lambda t - Q_v(t).\]
%Since $Q$ is of bounded variation by (\ref{31}), so are $Q_v$ and $L_v$. 
Define $\mathcal{R}_t$ and $\mathcal{Z}_t$ by the right hand sides of (\ref{r1}) and (\ref{r2}), respectively. We show that the pair of processes $(\mathcal{R}, \mathcal{Z})$ satisfies Zhang's fluid model. In fact, it suffices to show conditions (i) and (ii) of Zhang's fluid model.
%From the definition of $Q_v(t)$, we have for each $x\geq 0$, \begin{eqnarray} & &  \lambda (\bar{G}^r(x)-\bar{G}^r(x+ Q_v(t)/\la)) \non\\
%&=&  \int_{[0,(\chi(t)-t)^+]}\frac{g^r(y+t+x)}{\bar G^r(y)}\eta_0(dy)+ \la \int_{0}^{t\wedge \chi(t)} g^r(x+u) du. \label{Qvchire}
%\end{eqnarray}  

From the definition of $Q_v(t)$, we have for each $x\geq 0$,\[\la \int^{Q_v(t)/\la }_{t\wedge \chi(t)} \bar{G}^r(x+s) ds 
=  \bone_{\{\chi(t)\geq t\}}\int_{[0,\chi(t)-t]}\frac{\bar{G}^r(y+t+x)}{\bar G^r(y)}\eta_0(dy).\] 
Combining this, the construction of $\mathcal{R}$ and (\ref{f4}), we have 
\begin{eqnarray} \mathcal{R}_t(\RR_+) &=&  \bone_{\{\chi(t)\geq t\}}\int_{[0,\chi(t)-t]}\frac{\bar{G}^r(y+t)}{\bar G^r(y)}\eta_0(dy) + \la \int^{t\wedge \chi(t) }_{0} \bar{G}^r(s) ds \non\\ &=& \eta_t[0,\chi(t)] = Q(t). \label{Q*Q}
\end{eqnarray}  
Since $Q(t)$ is of bounded variation by (\ref{31}), the previous display implies that $Q_v$ and hence $L_v$ are also  of bounded variation. Thus, condition (i) of Zhang's fluid model holds.

Next we show that condition (ii) of Zhang's fluid model holds for $B^*(t)=\mathcal{Z}_t(\RR_+)$, $Q^*(t)=\mathcal{R}_t(\RR_+)$ and $X^*(t)=B^*(t)+Q^*(t)$. Note that for each $t\geq 0$, we have showed that $Q^*(t)=Q(t)$.  By using the definition of $B^*$, the construction of $\mathcal Z$ and the property of $\mathcal Z_0$, we have  \begin{eqnarray}B^*(t)&=& \mathcal{Z}_t(\RR_+) \nonumber\\ &=&   \mathcal{Z}_0(C_{t}) + \int_0^t \bar{G}^r(Q_v(s)/\la) \bar{G}^s(t-s) d L_v(s) \nonumber\\ &=& \int_{[0,H^s)}\frac{\bar G^s(y+t)}{\bar G^s(y)}\nu_0(dy) + \int_0^t \bar{G}^r(Q_v(s)/\la) \bar{G}^s(t-s) d L_v(s). \label{B*} \end{eqnarray}

%By using (\ref{Qvchire}) with $x=0$, we have
%\begin{eqnarray*} \lambda G^r(Q_v(t)/\la) & = &  \lambda (1-\bar{G}^r(Q_v(t)/\la)) \\
%&=&  \int_{[0,(\chi(t)-t)^+]}\frac{g^r(y+t)}{\bar G^r(y)}\eta_0(dy)+ \la \int_{0}^{t\wedge \chi(t)} g^r(u) du \\ &=& \int_{[0,\chi(t)]} h^r(x)\eta_t(dx). \end{eqnarray*}  

By \eqref{mrq}, (\ref{DisLimRtCond1}) and $\la(t) = \la$ for each $t\ge 0$, we have 
\begin{eqnarray*}
R(t) &=& 
\int_{[0,H^r)} \left(\int_{(y,y+t]} \bone_{\{y\leq \chi(x-y)-(x-y)\}}  dG^r(x)\right) \bar{G}^r(y)^{-1}\eta_0(dy)  \non\\
&&  \qquad + \la \int_0^t \int_{[0,H^r]} \bone_{\{x\leq s \wedge \chi(s)\}} dG^r(x) ds \non \\ &=& \int_{[0,H^r)} \left(\int_{(y,y+t]} \bone_{\{x\leq \chi(x-y)\}} dG^r(x)\right) \bar{G}^r(y)^{-1}\eta_0(dy)  \non\\
&&   \qquad 
+  \la \int_0^t G^r(\chi(s)\wedge s) ds \\
&=& \la \int_0^t G^r(Q_v(s)/\la)ds.
\end{eqnarray*}

In addition, since $G^r_d(Q_v(t)/\la) = Q(t)/\la$ by (\ref{Q*Q}), by the chain rule formula, 
$$Q(t) = Q(0)+ \int_0^t \bar G^r(Q_v(s)/\lambda)dQ_v(s).$$ 
These, together with (\ref{31}), imply that \[K(t)= \la t + Q(0)-Q(t)-R(t)=\int_0^t \bar{G}^r(Q_v(s)/\la) d L_v(s).\] Hence, by (\ref{B*}), $B^*(t)=B(t)$ and then $X^*(t)=X(t)$. Since $B, Q, X$ satisfy (\ref{f1})--(\ref{f12}), then $B^*,Q^*,X^*$ satisfy condition (ii) of Zhang's fluid model. This completes the proof that $(\mathcal{R}, \mathcal{Z})$ is a Zhang's fluid model. Clearly from the construction, both the given Zu{\~n}iga's fluid model and the constructed Zhang's fluid model describe the same $G/GI/N+GI$ queueing system. ~~~\bsq

\begin{coro} \label{coro:KRZ}
Given a Zu{\~n}iga's fluid model $(\eta,\nu,X)$ for a  $G/GI/N+GI$ queueing system with the input data $(\la,\eta_0,\nu_0,X(0))$ such that $(\eta_0,\nu_0,X(0))\in \mathcal S_0$ and $\eta_0(dx)= \la \bone_{[0,a]}(x) \bar{G}^r(x)dx$ for some $a\geq 0$, then one can construct a Zhang's fluid model $(\mathcal{R}, \mathcal{Z})$ for the same queueing system
with arrival rate $\la$, $\mathcal{Z}_0$ defined by \eqref{con2} for the given $\nu_0$, and $\mathcal{R}_0$ defined by \eqref{r3} for $Q_v(0)=a \la$. 
\end{coro}
\noindent {\it Proof.} It suffices to check that the given $\eta_0$ satisfies (\ref{DisLimRtCond}) and (\ref{DisLimRtCond1}). Note that for the given $\eta_0$, the equation in (\ref{DisLimRtCond}) becomes 
\begin{eqnarray*}
 \int^{z}_{t\wedge \chi(t)} \bar{G}^r(x+s) ds 
=  \bone_{\{\chi(t)\geq t\}}\int_t^{t+(\chi(t)-t)\wedge a}\bar{G}^r(s+x)ds.
\end{eqnarray*}
For $t\geq 0$ such that $\chi(t)\geq t$, we can choose $z_t=t+ (\chi(t)-t)\wedge a$ and for $t\geq 0$ such that $\chi(t)< t$,  we can choose $z_t=\chi(t)$. Clearly, in either case, $z_t$ does not depend on $x\geq 0$. Now we show that $z_t$ satisfies (\ref{DisLimRtCond1}). Note that for the given $\eta_0$, by (\ref{Rr1}),
\begin{eqnarray*}
& & \int_{[0,H^r)} \left(\int_{(y,y+t]} \bone_{\{x\leq \chi(x-y)\}} dG^r(x)\right) \bar{G}^r(y)^{-1}\eta_0(dy) \\ &=&  \int_{[0,H^r]} \int_{[0,t]} \bone_{\{x\leq \chi(s)\}} d_s \left(\int_{[(x-s)^+,x)} \frac{\bar{G}^r(x-)}{\bar{G}^r(y)} \eta_0(dy)\right) dM^r(x) \\ &=&  \la \int_{[0,H^r]} \int_{[0,t]} \bone_{\{x\leq \chi(s)\}} d_s \left(\int_{(x-s)^+}^x \bone_{[0,a]}(y) dy\right) dG^r(x) \\ &=&  \la \int_{[0,H^r]} \int_{[0,t]} \bone_{\{x\leq \chi(s)\}} \bone_{\{s<x\leq s+a\}} ds dG^r(x) \\ &=& \la \int_0^t (G^r(\chi(s)\wedge (s+a)) - G^r(s\wedge \chi(s))) ds.
\end{eqnarray*}
It follows that 
\begin{eqnarray*}
& & \la \int_0^t G^r(\chi(s)\wedge s) ds \non\\ & & \quad +  \int_{[0,H^r)} \left(\int_{(y,y+t]} \bone_{\{x\leq \chi(x-y)\}} dG^r(x)\right) \bar{G}^r(y)^{-1}\eta_0(dy) \\ &=& \la \int_0^t G^r(\chi(s)\wedge s) ds + \la \int_0^t (G^r(\chi(s)\wedge (s+a)) - G^r(s\wedge \chi(s))) ds \\ &=& \la \int_0^t \bone_{\{\chi(s)\geq s\}} \left(G^r(s)+G^r(\chi(s)\wedge (s+a)) - G^r(s) \right)ds \\ & & \quad + \la \int_0^t \bone_{\{\chi(s)< s\}} G^r(\chi(s))ds \\ &=& \la \int_0^t G^r(z_s)ds.
\end{eqnarray*}
 Thus, (\ref{DisLimRtCond1}) holds for the choice of $z_t$ and hence the corollary follows directly from Lemma \ref{lem:KRZ}. ~~~\bsq

\section{Proof of Theorem 2.1} \label{secConn}

In this section, we prove Theorem \ref{wellposed},  the equivalence between the two fluid models tracking elapsed times described in Sections \ref{secFluidmeas} and \ref{secFluidtwo}. 
We first derive a set of two-parameter fluid equations from a measure-valued fluid model $(\eta, \nu,  X)$ in Definition \ref{def:FMA} and show that it is a two-parameter fluid model; see Proposition \ref{fluidtwo}.
We then derive a set of measure-valued fluid equations from a two-parameter fluid model $(B(t,y),Q(t,y))$ in Definition \ref{def:FM2} and show that it is a measure-valued fluid model; see Proposition \ref{fluidtwo2}.
Thus we conclude that the existence and uniqueness of the two fluid models are equivalent.
% see Theorem \ref{wellposed}.
%As a consequence of this equivalence property, we show that the total fluid can be expressed as a differential equation, which reduces to the well-known ODE when the system is Markovian; see Proposition \ref{fluidrate} and Corollary \ref{fluidODE}. \footnote{not sure if they should be in this section or the section of main results...}

Recall that $\chi(t)$ in (\ref{f9}) represents the waiting time of the fluid content at the head of the queue. Namely, the fluid content in the potential queue must be in queue waiting for service if the waiting time is less than $\chi(t)$, but must have abandoned otherwise. 
By the FCFS service discipline, the definition of the potential queue and the role of $w(t)$, we see that 
 \beql{k1}
\chi(t)=w(t).
\eeq
 We also observe that evidently, for each $y\geq 0$,\beql{k2}
B(t,y)  = \nu_t[0,y], \quad \tilde{Q}(t,y)= \eta_t[0,y], \mbox{ and } Q(t,y)=\eta_t[0,y\wedge \chi(t)].
 \eeq

We first start with the measure-valued fluid model $(\eta, \nu,  X)$ in Definition \ref{def:FMA}, and show that the two-parameter processes $(B(t,y), Q(t,y))$ in \eqref{k2} satisfy  Definition \ref{def:FM2}. For this we need to assume that $\eta_0$ and $\nu_0$ have densities $\tilde q(0,x)$ and $b(0,x)$, respectively, since they are required in the definition of the two-parameter fluid model.

\begin{prop} \label{fluidtwo}
Let $(\eta,\nu,X)$ be a measure-valued fluid model tracking elapsed times with the input data $(\lambda(\cdot),\eta_0,\nu_0,X(0))$ such that $(\eta_0,\nu_0,X(0))\in \mathcal S_0$. 
Suppose that $\eta_0$ and $\nu_0$ have densities $\tilde q(0,x)$ and $b(0,x)$, respectively.  Then, $(B(t,y),Q(t,y))$ given by $(\ref{k2})$ is a two-parameter fluid model tracking elapsed times with the input data $(\lambda(\cdot), \tilde q(0,x), b(0,x))$ and $q(0,x) = \tilde q(0,x)$. 
\end{prop}
{\it Proof.} Let $(\eta,\nu,X)$ be a measure-valued fluid model tracking elapsed times with the input data $(\lambda(\cdot),\eta_0,\nu_0,X(0))$ such that $(\eta_0,\nu_0,X(0))\in \mathcal S_0$ and $\eta_0$ and $\nu_0$ have densities $\tilde q(0,x)$ and $b(0,x)$, respectively.  
Since $b(0,x)$ and $\tilde q(0,x)$ denote the densities of $\nu_0$ and $\eta_0$, respectively, it follows that $b(0,x) =0$ for each $x\geq H^s$ and  $\tilde q(0,x) =0$ for each $x\geq H^r$. For each $t\ge 0$ and $y \geq 0$, by letting $\psi_y(x) = \bone(0 \le x \le y)$   in \eqref{f4} and \eqref{f5}, respectively (Corollary 4.2 in \cite{KR10} shows that \eqref{f4} and \eqref{f5} hold for any bounded Borel measurable function $\psi$.),  $B(t,y)$ and $\tilde Q(t,y)$ satisfy the following equations, respectively:
\beqal{t2}
B(t,y) &=& \int_0^{(y-t)^+\wedge H^s} \frac{\bar{G}^s(x+t)}{ \bar{G}^s(x)} b(0,x) dx +\int_{(t-y)^+}^t \bar{G^s}(t-s) \kappa(s)ds \non\\ &=& \int_0^{(y-t)^+\wedge H^s} \frac{\bar{G}^s(x+t)}{ \bar{G}^s(x)} b(0,x) dx +\int_{0}^{y\wedge t} \bar{G^s}(s) \kappa(t-s)ds,
\eeqa

\beqal{t1}
\tilde{Q}(t,y) &=&  \int_0^{(y-t)^+\wedge H^r} \frac{\bar{G}^r(x+t)}{ \bar{G}^r(x)}\tilde q(0,x)dx +  \int_{(t-y)^+}^t \bar{G^r}(t-s) \la(s)ds  \non\\
&=&   \int_0^{(y-t)^+\wedge H^r} \frac{\bar{G}^r(x+t)}{ \bar{G}^r(x)} \tilde q(0,x)dx  + \int_0^{y\wedge t} \bar{G^r}(s) \la(t-s) ds.
\eeqa
Then from (\ref{t2}) and (\ref{t1}), $B(t,y)$ and $\tilde Q(t,y)$ have densities $b(t,y)$ and $\tilde q(t,y)$, respectively, with the representation:
\beql{bty}
b(t,y)=\left\{\begin{array}{ll} \bar{G^s}(y) \kappa(t-y) & \mbox{ if } y< t\wedge H^s, \\
 \frac{\bar{G}^s(y)}{ \bar{G}^s(y-t)} b(0,y-t)  & \mbox{ if } t<y<t+H^s, \\ 0 & \mbox{ otherwise, } \end{array}\right.
\eeq
and
\beql{qty}
\tilde q(t,y)=\left\{\begin{array}{ll} \bar{G^r}(y) \lambda(t-y) & \mbox{ if } y< t\wedge H^r, \\
 \frac{\bar{G}^r(y)}{ \bar{G}^r(y-t)} \tilde q(0,y-t)  & \mbox{ if } t<y<t+H^r, \\ 0 & \mbox{ otherwise. }  \end{array}\right.
\eeq
From this, it is easy to check that the two fundamental evolution equations in \eqref{w2} and \eqref{w3} are satisfied.
It is clear from the last equation in (\ref{k2}) that $Q(t,y)$ satisfies the following equation:

\beql{t3}
Q(t,y) =  \int_0^{(y\wedge \chi(t)-t)^+\wedge H^r} \frac{\bar{G}^r(x+t)}{ \bar{G}^r(x)} \tilde q(0,x)dx+\int_0^{y\wedge \chi(t) \wedge t} \bar{G^r}(s) \la(t-s) ds,
\eeq
Then, by comparing with (\ref{t1}), we have that
\beql{t6}
Q(t,y) = \left\{\begin{array}{ll}  \tilde{Q}(t,y) & \mbox{ if } y < \chi(t), \\
Q(t) & \mbox{ if } y \geq  \chi(t). \end{array}\right.
\eeq
Now, define $q(t,y)$ by
 \beql{tqty}
q(t,y) = \left\{\begin{array}{ll}  \tilde{q}(t,y) & \mbox{ if } y < \chi(t), \\
0 & \mbox{ if } y >  \chi(t),\\ \lambda(t)-\lambda(t)\wedge \int_{[0,H^s)}h^s(x)b(t,x) dx & \mbox{ if } y=\chi(t),\ B(t)=1,  \\ 0 & \mbox{ if } y=\chi(t), B(t)<1. \end{array}\right.
\eeq
Note that
\[\int_0^y q(t,x)dx=\int_0^y \tilde q(t,x)dx = \tilde Q(t,y)=Q(t,y), \mbox{ if } y<\chi(t), \] and \[\int_0^y q(t,x)dx=\int_0^{\chi(t)} \tilde q(t,x)dx = \tilde Q(t,\chi(t))=Q(t) =Q(t,y), \mbox{ if } y\geq \chi(t). \]
Thus, $q(t,y)$ is a density function of $Q(t,y)$. Since $(\eta_0,\nu_0,X(0))\in \mathcal S_0$, it is clear that $q(0,x)=\tilde q(t,x)$.

 From \eqref{f8} and \eqref{f9}, we obtain the following expression of $R(t)$  using the process $\tilde{Q}(t,y)$ in \eqref{t1},
\beqal{t4}
&& \qquad  R(t) = \int_0^t \left( \int_0^{\chi(s)\wedge H^r} h^r(x) \eta_s(dx) \right) ds \\
&=& \int_0^t \left(  \int_0^{(\chi(s)-s)^+ \wedge H^r}  \frac{g^r(x+s)}{\bar{G}^r(x)} \tilde q(0, x)dx+\int_0^{s\wedge\chi(s)} g^r(x)\la(s-x) dx \right) ds  \non \\ &=& \int_0^t \left( \int_0^{\chi(s) \wedge H^r} h^r(x) \tilde Q(t,dx)\right)\,ds, \non
\eeqa
and from (\ref{f80}), we obtain the following expression of $D(t)$ using the process $B(t,y)$ in (\ref{t2}),
\beqal{t70}
&& D(t) =  \int_0^t \left( \int_{[0,H^s)} h^s(x) \nu_s(dx) \right) ds \\
&=& \int_0^t \left(  \int_{[0,H^s)}  \frac{g^s(x+s)}{\bar{G}^s(x)} b(0, x)dx+\int_0^{s} g^s(x)\kappa(s-x) dx \right) ds  \non \\ &=& \int_0^t \left( \int_{[0,H^s)} h^s(x) B(t,dx)\right)\,ds. \non
\eeqa
From (\ref{t4}), (\ref{t70}) and (\ref{t6}), we can see that $D(t)$ and $R(t)$ have densities $\sigma(t)$ and $\alpha(t)$, respectively, and they satisfy \eqref{w6}, that is,
\beql{dr}
\sigma(t) = \int_{[0,H^s)} b(t,x) h^s(x) dx, \quad \alpha(t) = \int_{[0,H^r)} q(t,x) h^r(x) dx, \quad t\ge 0.
\eeq

To complete the proof, it is enough to show that $b(t,0)$, $\tilde q(t,0)$ and $q(t,0)$ from (\ref{bty}), (\ref{qty}) and (\ref{tqty}) satisfy (\ref{w7}), (\ref{w9}) and (\ref{w8}), respectively. Note that $b(t,0)=\kappa(t)$ by (\ref{bty}). Combining this with (\ref{kappa}) and (\ref{dr}), $b(t,0)$  satisfies (\ref{w7}). By (\ref{qty}), $\tilde q(t,0)=\lambda(t)$ and then satisfies (\ref{w9}). On the other hand,  from (\ref{tqty}) and \eqref{dr},
 \beql{qty0}
q(t,0) = \left\{\begin{array}{ll}  \tilde{q}(t,0)=\lambda(t) & \mbox{ if } 0 < \chi(t),  \\ \lambda(t)-\lambda(t)\wedge \sigma(t) & \mbox{ if } 0=\chi(t),\ B(t)=1, \\ 0 & \mbox{ if } 0=\chi(t),\ B(t)<1. \end{array}\right. \eeq  This implies that $q(t,0)$ satisfies (\ref{w8}). Finally, the rate balance equation \eqref{w10} follows from the balance equation \eqref{31}, by noting that $Q(t) = \int_0^{\infty} q(t, x) dx$, $K(t) = \int_0^t b(s,0) ds$ and $R(t) = \int_0^t \alpha(s) ds$.   ~~~\bsq

\bigskip

We next show that a set of measure-valued equations $(\nu, \eta, X)$ derived from a two-parameter fluid model in Definition \ref{def:FM2} satisfies Definition \ref{def:FMA}.

\begin{prop} \label{fluidtwo2}
Let $(B(t,y),Q(t,y))$ be a two-parameter fluid model tracking elapsed times with the input data $(\lambda(\cdot), \tilde q(0,x),b(0,x)) $ and $q(0,x) =  \tilde q(0,x) $. For each $t\geq 0$, let $\eta_t[0,y]\doteq \tilde Q(t,y)$ and $\nu_t[0,y] \doteq B(t,y)$ for each $y\geq 0$ and define $X(t)\doteq B(t,\infty)+Q(t,\infty)$. Then, $(\eta,\nu,X)$ is a measure-valued  fluid model tracking elapsed times with the input data $(\lambda(\cdot),\eta_0,\nu_0,X(0))$ such that $(\eta_0,\nu_0,X(0))\in \mathcal S_0$. \end{prop}
\textit{Proof}. Fix $(B(t,y),Q(t,y))$ and the triple of functions $(\eta, \nu, X)$ defined from it. It is clear from the two fundamental evolution equations (\ref{w2}) and (\ref{w3}) that for each $t\geq 0$, $\tilde q(t,x)$ as a function in $x$ has support in $[0,H^r)$ and $b(t,x)$ as a function in $x$ has support in $[0,H^s)$. It then follows that $\eta_t$ has support in $[0,H^r)$ and $\nu_t$ has support in $[0,H^s)$ for each $t\geq 0$.  Also it is clear that $(\eta_0,\nu_0,X(0))\in \mathcal S_0$.

 We first show that $\nu$ satisfies (\ref{f5}).  For every $\psi \in \sC_b(\RR_+)$  and $t\geq 0$,
\beql{p1}
\int_0^{\infty} \psi(x) \nu_t (d x) =  \int_0^{\infty} \psi(x)b(t,x) dx =  \int_0^{t} \psi(x)b(t,x) dx + \int_t^{\infty} \psi(x)b(t,x) dx.
\eeq
For the first term on the right-hand side of (\ref{p1}), we can use the first fundamental evolution equation (\ref{w2}) to yield that
\beql{p2}
\int_0^{t} \psi(x)b(t,x) dx =  \int_0^{t} \psi(x)b(t-x,0)\frac{\bar G^s(x)}{\bar G^s(0)} dx =  \int_0^{t} \psi(x)\bar G^s(x)b(t-x,0)dx .
\eeq
For the second term on the right-hand side of (\ref{p1}), another application of the first fundamental evolution equation (\ref{w2}) yields that
\beqal{p3}
\int_t^{\infty} \psi(x)b(t,x) dx &=&  \int_{t\wedge H^s}^{H^s} \psi(x)b(t,x) dx \\ &=&  \int_{t\wedge H^s}^{H^s} \psi(x)b(0,x-t)\frac{\bar G^s(x)}{\bar G^s(x-t)} dx \non \\
&=&  \int_0^{t\vee H^s-t} \psi(x+t)\frac{\bar G^s(x+t)}{\bar G^s(x)}b(0,x)dx \non \\ &=&  \int_0^{H^s} \psi(x+t)\frac{\bar G^s(x+t)}{\bar G^s(x)}b(0,x)dx. \non 
\eeqa
The last equality in (\ref{p3}) follows from the fact that $\bar G^s(x+t)=0$ if $x\in (t\vee H^s-t,H^s)$. 
For each $t\geq 0$, let \beql{p5}
K(t)\doteq \int_0^tb(s,0)ds.
\eeq Combining the above four displays, we obtain that $\nu$ satisfies (\ref{f5}). An analogous argument using the second fundamental evolution equation (\ref{w3}) shows that
\beqal{p4}
&& \int_0^{\infty} \psi(x) \eta_t (d x) =  \int_0^{\infty} \psi(x)\tilde q(t,x) dx  \\ &=&   \int_0^{t} \psi(x)\tilde q(t,x) dx + \int_t^{\infty} \psi(x)\tilde q(t,x) dx \non \\ &=& \int_0^{t} \psi(x)\bar G^r(x)\tilde q(t-x,0)dx+ \int_0^{H^r} \psi(x+t)\frac{\bar G^r(x+t)}{\bar G^r(x)}\tilde q(0,x)dx. \non
\eeqa
By (\ref{w9}), $\tilde q(t-x,0)=\lambda(t-x)$. Thus, $\eta$ satisfies (\ref{f4}).

Next, for each $t\geq 0$, define $B(t)\doteq B(t,\infty)$, $Q(t)\doteq Q(t,\infty)$, $D(t)\doteq \int_0^t\sigma(x)dx$, $R(t)\doteq \int_0^t \alpha(x)dx$. From (\ref{w6}), $D$ satisfies (\ref{f80}).
 From
(\ref{w6}) again, (\ref{w5}) and (\ref{k1}), $R$ satisfies (\ref{f8}). Since $\nu$ satisfies (\ref{f5}), by choosing $\psi=\bone$ in (\ref{f5}), we have \[B(t)= \int_0^{H^s}\frac{\bar{G^s}(x+t)}{\bar{G^s}(x)}\nu_0(dx)+
 \int_0^t  \bar{G^s}(t-s) d K(s)\]
and by choosing $\psi=h^s$ in (\ref{f5}), we have \begin{eqnarray*}
D(t)&=& \int_0^t \int_0^{H^s} b(s,x)h^s(x)dx ds \\ &=&  \int_0^t \left( \int_0^{H^s}\frac{g^s(x+s)}{\bar{G^s}(x)}\nu_0(dx)+
 \int_0^s  g^s(s-x) d K(x)\right)ds \\ &=&   \int_0^{H^s}\frac{\bar{G^s}(x)-\bar{G^s}(x+t)}{\bar{G^s}(x)}\nu_0(dx)+\int_0^t
 \int_0^s  g^s(s-x) d K(x)ds \\ &=& B(0)- \int_0^{H^s}\frac{\bar{G^s}(x+t)}{\bar{G^s}(x)}\nu_0(dx)+ \int_0^tG^s(t-s)dK(s) \\ &=& B(0)-B(t)+K(t).
\end{eqnarray*}
This shows that (\ref{33}) is satisfied. The relationship (\ref{w10}) directly implies that (\ref{31}).

Notice from (\ref{w8}) that $q(t,0)=0$ when $B(t)<1$ and $q(t,0)=\lambda(t)>0$ when $Q(t)>0$. Thus, $B(t)<1$ and $Q(t)>0$ cannot happen at the same time since $\lambda(t)>0$. As the consequence, (\ref{f12}) holds. This also shows that \[(X(t)-1)^+=(B(t)+Q(t)-1)^+=Q(t),\] 
and
\[ B(t)=X(t)\wedge 1= 1-(1-X(t))^+,\] 
which imply that both (\ref{f1}) and (\ref{f2}) hold. At last, (\ref{f9}) follows directly from (\ref{k1}). This completes the proof of the proposition. ~~~\bsq
\bigskip

\section*{Acknowledgement}
%We thank the referees for their helpful comments and suggestions, which have helped the quality and exposition of the paper.  
We thank Ward Whitt and Kavita Ramanan for many helpful discussions with the paper which have helped improve the exposition of the paper.

\singlespacing

%\section{Appendix}
\appendix

\renewcommand{\theequation}{B.\arabic{equation}}
\renewcommand\thetheorem{B.\arabic{theorem}}

 \end{document}